\documentclass[a4paper]{article}
\usepackage{etoolbox}
\usepackage{abstract}
\newbool{PREPRINT}
\booltrue{PREPRINT}
\usepackage{makecell}
\usepackage{soul}
\usepackage[colorlinks,bookmarksopen,bookmarksnumbered,citecolor=red,urlcolor=red]{hyperref}
\usepackage{stmaryrd}
\usepackage[usenames]{color}
\usepackage{comment}
\usepackage{amsmath,amssymb,amsthm}
\usepackage{mathrsfs}
\usepackage{graphicx}
\usepackage{algorithm}
\usepackage{algorithmicx}
\usepackage{algcompatible}
\usepackage{caption}
\usepackage{listings}
\usepackage{mathtools}
\usepackage{multicol}
\usepackage{placeins}
\usepackage{multirow}
\usepackage{subcaption}
\usepackage{enumitem}
\usepackage{fdsymbol}
\setlist[description]{itemsep=-0.5ex}
\newcommand{\figdir}{./}

\newif\ifpreprint
\ifbool{PREPRINT}{ 
\preprinttrue
\usepackage[cm]{fullpage}
\usepackage{amsthm}
\usepackage{multicol}
\setlength\columnsep{20pt}
\usepackage[affil-it]{authblk}

}{ %
\usepackage{fullpage}
\newdefinition{defn}{Definition}
\newdefinition{expl}{Example}
} 

\newcommand{\lag}{\mathcal{L}}
\newcommand{\ham}{\mathcal{H}}

\newtheorem*{theorem}{Theorem}

\lstdefinelanguage[ppmd]{python}[]{python}{%
  emph={ParticleLoop,ParticleDat,PositionDat,ScalarArray,GlobalArray,Kernel,PairLoop,Constant,State,Data,IntegratorRange}
}
\definecolor{DarkBlue}{rgb}{0.00,0.00,0.55}
\definecolor{DarkRed}{rgb}{0.55,0.00,0.00}
\definecolor{DarkGreen}{rgb}{0.00,0.55,0.00}
\definecolor{Gray}{rgb}{0.95,0.95,0.95}
\definecolor{Purple}{rgb}{0.5,0.0,0.5}
\definecolor{Bittersweet}{rgb}{1.0,0.44,0.37}
\newcommand{\secapp}{\ifbool{PREPRINT}{Appendix~}{}}
\newcommand{\pparagraph}[1]{\ifbool{PREPRINT}{\paragraph{#1.}}{\paragraph{#1}}}

\title{Exact conservation laws for neural network integrators of dynamical systems}
\ifbool{PREPRINT}{ 
  \author[a,*]{Eike~Hermann~M\"{u}ller}
  \affil[a]{University of Bath, Bath BA2 7AY, Bath, United Kingdom}
\affil[*]{Email: \texttt{e.mueller@bath.ac.uk}}
}{
\author[1]{Eike~Hermann~M\"{u}ller\fnref{fn1}}
\ead{e.mueller@bath.ac.uk}
\fntext[fn1]{email: \texttt{e.mueller@bath.ac.uk}}
\address[1]{Department of Mathematical Sciences, University of Bath, Bath BA2 7AY, Bath, United Kingdom}
} 

\begin{document}
\ifbool{PREPRINT}{ 
  \maketitle
  \begin{onecolabstract}
}{%
  \begin{abstract}
    } 
    The solution of time dependent differential equations with neural networks has attracted a lot of attention recently. The central idea is to learn the laws that govern the evolution of the solution from data, which might be polluted with random noise. However, in contrast to other machine learning applications, usually a lot is known about the system at hand. For example, for many dynamical systems physical quantities such as energy or (angular) momentum are exactly conserved. Hence, the neural network has to learn these conservation laws from data and they will only be satisfied approximately due to finite training time and random noise. In this paper we present an alternative approach which uses Noether's Theorem to inherently incorporate conservation laws into the architecture of the neural network. We demonstrate that this leads to better predictions for three model systems: the motion of a non-relativistic particle in a three-dimensional Newtonian gravitational potential, the motion of a massive relativistic particle in the Schwarzschild metric and a system of two interacting particles in four dimensions.
    \ifbool{PREPRINT}{ 
    \end{onecolabstract}
    \textbf{keywords}:
    \newcommand{\sep}{, }
    }{
  \end{abstract}
  \begin{keyword}
    } 
    Machine Learning\sep Lagrangian Mechanics\sep Noether's Theorem\sep Dynamical System
    \ifbool{PREPRINT}{ 
    \\[1ex]
    }{
  \end{keyword}
  \maketitle
}
\section{Introduction}
Traditionally, machine learning has been applied to areas where there are no known underlying laws that describe the observed data. For example, deep convolutional neural networks (CNNs) have been very successful in computer vision applications (see \cite{He2016,Ren2015,Krizhevsky2012} for some famous examples), even though it is unclear how the lower-dimensional manifold of ``valid'' images is parametrised. This manifold has to be discovered in a purely data-driven way by feeding the network vast amounts of data. Of course, also in this case domain specific prior knowledge is built into the structure of the network: deep CNNs are so successful since the network structure supports the hierarchical representation of higher level features and encapsulates the fact that the input images are not random, but nearby pixels are likely to have similar values. In contrast, other fields of science such as physics are governed by well known fundamental laws. There appears to be little value in applying machine learning in this case since these laws allow the construction of an exact mathematical description, which can then be translated into an algorithm for making predictions with traditional techniques from numerical analysis. However, often the dynamics of a system is only partially constrained by these underlying laws. For example, the motion of a particle in an external potential might conserve energy and angular momentum, while the exact form of the potential or the expression for the kinetic energy is still unknown and needs to be inferred from observations. In fact, a very successful approach to constructing new theories in modern physics is to write down the most general Lagrangian which is invariant under certain symmetry transformations and to then constrain the remaining parameters through a fit to experimental data. Famously, general relativity \cite{Einstein1923} is constructed by demanding that the theory is invariant under local coordinate transformations but this allows the existence of a cosmological constant which needs to be constrained by experimental evidence. In particle physics a similar method is used for the construction of Chiral Perturbation Theory \cite{Gasser1984,Scherer2003} and other effective theories for which the Lagrangian can be written as an infinite sum of terms that are invariant under certain symmetry transformations. The expansion coefficients have to be inferred from experiments.

In this work we study simpler systems arising in classical mechanics of point particles. The key idea is to represent the Lagrangian by a neural network in such a way that it exactly conserves quantities like linear and angular momentum. The weights of the network are then learned from training data which consists of trajectories that follow the true dynamics and are perturbed by random noise. While in the present paper we use synthetic training data that is obtained by integrating the exact equations of motion numerically for a set of model systems, it should be stressed that the data could also come from the observation of a real physical system, which naturally includes measurement errors. In our numerical experiments these errors are modeled by adding a small normally distributed random variable to the exact trajectory. Drawing a historical analogy, our neural network can be seen as the virtual version of an astronomer who tries to derive the celestial equations of motion from (inexact) observations of the planets in the night sky. In contrast to its medieval predecessors, our virtual astronomer has read Noether's paper \cite{Noether1918} and takes care not to violate fundamental physical laws such as angular momentum conservation. In contrast to \cite{Pfaff2020,Sanchez2020} where noise is used to regularise the training process, here we assume that for real-life systems noise is inherent to the measurement process and in principle outside the control of the observer.
\subsection{Related work}\label{sec:relatedwork}
Several approaches for solving time dependent differential equations with neural networks have been pursued in the literature. In \cite{Kadupitiya2022} the authors train an LSTM \cite{Hochreiter1997} based neural network to predict the solution at the next time step, given the solution at a number of $s$ previous times. Hence, this can be seen as an extension of classical linear multistep methods (see e.g. \cite{Hairer1993}) and the hope is that with a suitably trained network it is possible to substantially increase the timestep size without losing accuracy, thereby making the method more efficient than traditional methods.
\subsubsection{Hamiltonian- and Lagrangian neural networks}
For Hamiltonian systems a completely different approach is pursued in \cite{Greydanus2019}: a neural network is trained to learn a scalar valued function $\ham_{\text{NN}}(q,p)$ which approximates the true Hamiltonian $\ham(q,p)$ as a function of the generalised coordinates $q\in\mathbb{R}^{d}$ and conjugate momenta $p\in\mathbb{R}^{d}$; here and in the following $d$ is the dimension of the dynamical system. It is well known that symplectic integrators preserve a so-called shadow Hamiltonian $\widetilde{\ham}$, i.e. the numerical solution generated with the true Hamiltonian $H$ is the \textit{exact} solution of $\widetilde{\ham}$. For this reason, the numerical solution approximately preserves energy. While the authors on \cite{Greydanus2019} do not exploit this since they use a non-symplectic fourth order Runge Kutta integrator, in \cite{Chen2019} it is argued that symplectic Neural Networks (SRNNs) show superior performance. SRNNs represent both the kinetic and potential energy by a neural network and then use a symplectic Verlet method \cite{Verlet1967} to propagate the solution. Instead of working in the Hamiltonian formulation, it might be more convenient to represent the Lagrangian by a neural network as in \cite{Cranmer2020}. As argued there, working in the Lagrangian framework is more flexible since the system can be formulated in \textit{any} coordinate system, not just canonical coordinates, which are often not easy to find. For medium to high-dimensional systems one issue is that the Lagriangian depends on many variables and errors are amplified when taking derivatives, which leads to inaccurate predictions for the accelerations. Furthermore, the inversion of the Hessian matrix might not only become numerically unstable but also prohibitively expensive.
\subsubsection{Physics-motivated networks with built-in symmetries}
The methods discussed so far only preserve some physical quantities such as the total energy. It can be shown that energy conservation directly follows from time invariance of the system, i.e. the fact that the Lagrangian does not explicitly depend on time. Noether's Theorem \cite{Noether1918} generalises this result and shows that invariance of the Lagrangian under \textit{any} infinitesimal transformation results in a corresponding conservation law. Hence, if the neural network that represents the system's Lagrangian can be constructed such that it is exactly invariant under a certain infinitesimal symmetry transformation, this will guarantee the conservation of a physical quantity. The role of symmetries in neural networks has been explored in \cite{Mattheakis2019}. The authors use so-called ``hub-neurons'' to construct a neural network which represents a symmetric function $f:\mathbb{R}\rightarrow \mathbb{R}$ that satisfies $f(x)=f(-x)$. They then show that this gives a better fit to ground-truth training data which is polluted by non-symmetric noise. They use a similar technique to construct a symplectic neural network for the solution of time-dependent systems, but it is not clear whether their network is truly symplectic. The central role of symmetries is also exploited in \cite{Ling2016}. Here the aim is to use neural networks to represent closures for the Reynolds-averaged Navier Stokes equations. Such a closure allows computing the shear-stress in the fluid from the non-turbulent component of the velocity field. Crucially, if the flow is isotropic, the shear stress can only depend on combinations of the velocity field which are invariant under rotations. The authors of \cite{Ling2016} achieve this by systematically constructing rotationally invariant scalars from the vector field and restricting the input of their neural network to these scalars.

The construction of neural networks that conserve physical quantities has been explored in other contexts. In \cite{Beucler2021} the authors outline a general procedure for choosing dynamical variables in such a way that they observe constraint equations which are conserved exactly  by construction; the method is applied to the parametrisation of convection and cloud physics in a climate model.
The authors of \cite{Sturm2020,Sturm2022} describe how mass conservation can be enforced in photochemical atmospheric models, which involve the transformation between different chemical species. Instead of predicting the change in concentration directly, the unknowns in the pen-ultimate layer of the neural network in \cite{Sturm2022} represent the flux-integrals, which are then combined into concentration changes with a known reaction. The authors show that this leads to exact mass conservation for some atomic species and improves the quality of prediction in general, both by reducing the error and by minimising unphysical negative outputs.
\subsubsection{Construction of equivariant and invariant functions}
The role of equivariance in the design of functions that model physical systems is explored in detail in \cite{Villar2021}. Building on classical work in \cite{Weyl1946}, the authors describe the systematic construction of functions that are equivariant under the action of Lie groups such as rotations or Lorentz- and Poincar\'e- transformations. Invariance can be seen as a special case of equivariance where the output transforms under the trivial representation of the symmetry group. The symmetry-enforcing layers that we construct in this paper are invariant functions in the language of \cite{Villar2021}. While in \cite{Villar2021} equivariant models are used to extract physical quantities from static data, the approach is extended to dynamical systems in \cite{Yao2021} where both the forcing function and the Hamiltonian of a springy pendulum system are constructed in an equivariant way and the authors find that this improves predictions. However, in constrast to our work, \cite{Yao2021} does not discuss the impact of invariance on conserved quantities.
\subsubsection{Equivariant convolutional neural networks}
Another approach to constructing rotation- and translation-equivariant neural networks is described in \cite{Thomas2018}. Following \cite{Cohen2016}, the authors focus on \textit{convolutional} layers and generalise the concept of translation-invariant convolutions (which are crucial for the success of CNNs in image processing) to rotation groups. The central idea is to write the filters as a product of a radial part and a spherical harmonic which transforms under a known representation of the group and to decompose the output into irreducible representations using Clebsch-Gordon coefficients. Compared to \cite{Villar2021}, which attempts to write down the most general equivariant function, the layers in \cite{Thomas2018} have a convolutional structure. In \cite{Thomas2018} the method is applied to extract the moment of inertia tensor (which transforms according to a 5-dimensional irreducible representation) from a point-mass distribution. Equivariant convolutional layers are also employed in \cite{Finzi2020}, which avoids the use of representation theory and incorporates locality into the network. The method is applied to simulate the of trajectories of six interacting particles. The authors explicitly discuss the conservation of linear and angular momentum with reference to Noether's Theorem and demonstrate that the invariance of the network leads to much better conservation of momentum (see \cite[Appendix A.5]{Finzi2020}). In contrast to our work, they use the Hamiltonian formulation in \cite{Greydanus2019}. They also impose the (strong) constraints that the Hamiltonian is separable and that the kinetic energy has the form $\sum_{j=1}^{N}p_j^2/(2m_j)$, where $p_j$ and $m_j$ are the particle momenta and masses respectively. We do not make any such assumptions here.
\subsubsection{Graph Neural Networks}
Graph Neural Networks (GNNs) \cite{Battaglia2018,Battaglia2016} attempt to incorporate inductive physical biases through intermediate layers that are modeled on a connectivity graph. In each layer information is propagated between the vertices (= object locations) and edges (= object relationships) of a graph network, which is motivated by the relationship between entities of the simulation. GNNs are therefore particularly well suited for modelling systems of interacting objects or particles. Embeddings are used to map physical quantities to the nodes and vertices of the network in the first layer. The GNN architecture allows the introduction of further constraints, for example locality in emulators for smooth particle hydrodynamics \cite{Sanchez2020}. Furthermore, compared to naive dense network layers, the number of parameters is dramatically reduced for systems of many particles since the activation functions share parameters across all nodes and edges. This advantage of GNNs for large systems is demonstrated in \cite{Sanchez2020} which considers emulators for thousands of particles. However, in contrast to our approach, the network in \cite{Sanchez2020} models the forcing functions directly (instead of returning the Lagrangian, from which the forces can be obtained via differentiation). By choosing suitable edge embeddings in the ``relative'' encoder in \cite{Sanchez2020}, the GNN attempts to capture translational invariance, which indeed improves prediction accuracy. However, as argued in \cite[Appendix A]{Finzi2020}, equivariance of the forcing function does not necessarily imply the existence of a conserved quantity. This can only be guaranteed if the Hamiltonian (or equivalently the Lagrangian) is invariant under a continuous symmetry, and hence the dynamics of the particles in \cite{Sanchez2020} is not characterised by a quantity that is conserved exactly, which is a key difference to our work. In \cite{Sanchez2019} the GNN architecture is combined with the Hamiltonian neural network approach \cite{Greydanus2019}, and it is demonstrated that this can lead to more accurate predictions. Although the authors point out that - assuming a symplectic integrator is used for timestepping - this will lead to the (approximate) conservation of energy, the conservation of other quantities (such as linear or angular momentum) is not discussed there. The focus of our work is to explore the relationship between invariant neural networks and conservation laws and we have used simple and general dense layers. Our approach will therefore not necessarily scale to systems of many interacting particles, similar to how imaging applications require CNNs instead of dense networks. We believe that a simple, dense neural network architecture is appropriate for the relatively small problems considered here. Inspired by \cite{Sanchez2019}, we outline some ideas on how GNNs could be combined with the approach discussed in this paper to simulate larger systems in Section \ref{sec:conclusion}.
\subsubsection{Neural networks based on the GENERIC formalism}
In \cite{Zhang2022} GENERIC Formalism Informed Neural Networks (GFINNs) are used to simulate systems that include both conservative and dissipative processes. The GENERIC approach \cite{Grmela1997,Oettinger1997,Oettinger2005,Oettinger2018} models these systems by two functions $E$ and $S$ (for the energy and entropy respectively) together with two matrices $L$ and $M$. The first and second laws of thermodynamics ($dE/dt=0$ and $dS/dt>0$) are satified if $E$, $S$, $L$ and $M$ obey certain conditions. In \cite{Zhang2022} these conditions are embedded into the architecture of the neural networks that are used to model the functions $E$, $S$ and the matrices $L$, $M$. As a result, the predicted trajectories are physical in the sense that they exactly preserve energy and the entropy is non-decreasing. The performance of GFINNs is also compared to two similar GENERIC based approaches \cite{Hernandez2021,Lee2021} for a range of model systems. Structure Preserving Neural Networks (SPNNs) \cite{Hernandez2021} assume that $E$ and $S$ are known and only model $L$ and $M$; the degeneracy conditions $M\,dE/dz=0$ and $L\,dE/dz=0$ are only enforced weakly through a term in the loss function which means that the first two laws of thermodynamics do not hold exactly. In contrast, GENERIC Neural Ordinary Equations (GNODEs) \cite{Lee2021} use a parametrisation of the bracket structure that enforces the degeneracy conditions exactly. While compared to the work presented here GENERIC based approaches \cite{Zhang2022,Hernandez2021,Lee2021} consider more general problems\footnote{In the GENERIC formalism Hamiltonian systems correspond to the special case $L=\begin{pmatrix}0&1\\-1&0\end{pmatrix}$, $M=0$ \cite[Remark 3.1]{Lee2021}}, they do not discuss the conservation of other physical quantities such and linear or angular momentum.
\subsection{Aim of this paper}\label{sec:achievements}
In this paper we show how neural networks can be used to solve Lagrangian dynamical systems while exactly preserving physical quantities. In contrast to \cite{Mattheakis2019} we consider invariance under \textit{continuous} transformations. For this, we follow \cite{Cranmer2020} and represent the Lagrangian as a neural network $\lag_{\text{NN}}(q,\dot{q})$ which takes as input the coordinates $q$ and velocities $\dot{q}$. As in \cite{Ling2016}, the first layer of the network reduces $q,\dot{q}$ to a set of scalar variables which do not change under the symmetry transformations that correspond to the conserved quantities. Noether's Theorem then guarantees exact conservation and during training the network learns the finer details of dynamics of the system that are \textit{not} constrained by the conservation laws. Since the network has to find solutions in a smaller sub-manifold of the entire solution space, training is likely to be more efficient. The main achievements of this paper are as follows:
\begin{enumerate}
  \item We show how invariance under continuous symmetries can be built into Lagrangian neural network models by passing the input through a so-called symmetry-enforcing layer (which can be seen as a collection of invariant functions \cite{Villar2021}); according to Noether's Theorem this leads to the exact conservation of a corresponding physical quantity.
  \item To demonstrate this, we consider three example systems:
        \begin{enumerate}
          \item The motion of a single particle in three dimensions under the influence of a central force field, which is invariant under rotations and thus conserves angular momentum.
          \item The motion of a massive particle in the rotationally invariant Schwarzschild metric, which is the relativistic pendant of the first problem. Again, three dimensional angular momentum is conserved.
          \item The motion of two interacting particles in $D$-dimensional space, where the interaction potential is invariant under rotations and translations. This results in the conservation of the $D$-dimensional linear momentum vector and the $\frac{1}{2}D(D-1)$ independent components of the antisymmetric $D\times D$ angular momentum tensor.
        \end{enumerate}
  \item For all systems we demonstrate that the generated trajectories are more realistic and that they conserve angular momentum and (in the case of the two-particle system) linear momentum to a high degree of accuracy.
  \item We further show numerically that enforcing invariance under continuous symmetry transformations makes the trajectories more stable with respect to small perturbations of the initial conditions.
\end{enumerate}
\paragraph{Structure} This paper is organised as follows: in Section \ref{sec:method} we review Noether's Theorem and explain the construction of Lagrangian neural networks that are invariant under continuous symmetry transformations. We write down the explicit form of the symmetry-enforcing input layers for the three model systems considered in this work. Numerical experiments for the model systems are described in Section \ref{sec:results}, where we present results that demonstrate the superior performance of our approach. Section \ref{sec:conclusion} contains our conclusions and ideas for future work. Some more technical details are relegated to the appendices: we discuss the construction of scalar invariants under the special orthogonal group $SO(D)$ in \ref{sec:rotationally_invariant_combinations} and present loss histories for the different neural networks in \ref{sec:loss_histories}.
\section{Methodology}\label{sec:method}
\subsection{Continuous symmetries and conservation laws}\label{sec:noether_theorem}
For completeness and further reference, we start by writing down Noether's Theorem here and refer the reader to \cite{Arnold2013} for a proof and further details.
\begin{theorem}[Noether \cite{Noether1918}]
  Consider a dynamical system which is formulated on a manifold $M$ with tangent bundle $TM$; the Lagrangian of this system is a real-valued function $\lag:TM\rightarrow \mathbb{R}$. Let $h^s:M\rightarrow M$ be a family of maps parametrised with the continuous parameter $s\in\mathbb{R}$ and let $h^s_{*,q}:TM_q\rightarrow TM_{h^s(q)}$ be the local derivative which maps between the tangent bundles at $q$ and $h^s(q)$. If the Lagrangian $\lag$ is invariant under $h^s$, i.e.
  \begin{equation}
    \lag(h^s(q),h^s_{*,q}(\dot{q})) = \lag(q,\dot{q})\qquad\text{for all $q,\dot{q}\in TM$},\label{eqn:L_invariance}
  \end{equation}
  then the following quantity is a constant of motion in the sense that $dI/dt=0$:
  \begin{equation}
    I = \frac{\partial \lag}{\partial \dot{q}} \frac{dh^s(q)}{ds}\Big|_{s=0}.\label{eqn:noether_conserved_I}
  \end{equation}
\end{theorem}

As an example, consider a particle moving in three dimensional space $\mathbb{R}^3$ and assume that the Lagrangian is invariant under rotations, i.e. the action of the special orthogonal group $SO(3)$. In this case, $M=\mathbb{R}^3$, $TM=\mathbb{R}^3\times \mathbb{R}^3$ and the maps $h^{s}$, $h^{s}_{*,q}$ can be written as linear transformations
\begin{xalignat}{2}
  h^{s}(q) &= \exp\left[s \Gamma \right]q, & h^{s}_{*,q}(v) &= \exp\left[s \Gamma \right]v
  \qquad\text{for $q,v\in\mathbb{R}^3\times\mathbb{R}^3$}
  \label{eqn:rotations}
\end{xalignat}
where $\Gamma$ is an element of the fundamental representation of the Lie-algebra $\mathfrak{so}(3)$ of $SO(3)$. The three matrices corresponding to rotations around the $x$-, $y$- and $z$- axis are
\begin{xalignat}{3}
  \Gamma_x &= \begin{pmatrix}0 & 0 & 0\\0 & 0 & -1\\0 & 1 & 0\\\end{pmatrix}, &
  \Gamma_y &= \begin{pmatrix}0 & 0 & 1\\0 & 0 & 0\\-1 & 0 & 0\\\end{pmatrix}, &
  \Gamma_z &= \begin{pmatrix}0 & -1 & 0\\1 & 0 & 0\\0 & 0 & 0\\\end{pmatrix}.
\end{xalignat}
Further, we have that
\begin{equation}
  \frac{dh^{s}}{ds}\Big|_{s=0} = \Gamma q
\end{equation}
and hence there are three conserved quantities, namely
\begin{xalignat}{3}
  L_x &= \frac{\partial \lag}{\partial \dot{q}_3}q_2 - \frac{\partial \lag}{\partial \dot{q}_2}q_3, &
  L_y &= \frac{\partial \lag}{\partial \dot{q}_1}q_3 - \frac{\partial \lag}{\partial \dot{q}_3}q_1, &
  L_z &= \frac{\partial \lag}{\partial \dot{q}_2}q_1 - \frac{\partial \lag}{\partial \dot{q}_1}q_2.
  \label{eqn:conserved_L_rotational}
\end{xalignat}
For a non-relativistic particle such as the one discussed in Section \ref{sec:kepler_problem}, the kinetic energy in the Lagrangian is $\frac{m}{2}(\dot{q}_1^2+\dot{q}_2^2+\dot{q}_3^2)$ and hence $\partial L/\partial \dot{q}_j = m\dot{q}_j$. In this case, the three conserved quantities are the components of the angular momentum vector
\begin{equation}
  L = m\;q \times \dot{q}.\label{eqn:angular_momentum_3d}
\end{equation}
It is worth pointing out that the quantities in Eq. \eqref{eqn:conserved_L_rotational} are conserved for more general cases, such as for the Lagrangian that governs the motion of a massive relativistic particle in the rotationally invariant Schwarzschild metric discussed in Section \ref{sec:schwarzschild_problem}.
\subsection{Neural networks with built-in conservation laws}\label{sec:nn_with_symmetries}
The central idea of this paper is as follows: we wish to construct a mapping $\Phi:\mathbb{R}^{2d}\rightarrow \mathbb{R}^d$ which predicts the acceleration $\ddot{q}=\Phi(q,\dot{q})$ from the position $q$ and velocity $\dot{q}$ such that the dynamics of this system preserves certain quantities \textit{exactly}. To achieve this, we use the Lagrangian formalism and represent the Lagrangian $\lag_{\text{NN}}$ by a neural network as in \cite{Greydanus2019}. The equations of motion are obtained by finding the stationary points of the action (= the time integral of the Lagrangian) which implies that
\begin{equation}
  \frac{d}{dt}\frac{\partial \lag_{\text{NN}}}{\partial\dot{q}} - \frac{\partial \lag_{\text{NN}}}{\partial q} = 0.\label{eqn:lagrangian_stationary}
\end{equation}
Following \cite[Eq. (6)]{Cranmer2020} the acceleration can then be obtained by taking the total time derivative of Eq. \eqref{eqn:lagrangian_stationary}:
\begin{equation}
  \ddot{q} = J_{\dot{q},\dot{q}}^{-T}\left( \frac{\partial \lag_{\text{NN}}}{\partial q} - J^T_{q,\dot{q}}\dot{q}\right) \label{eqn:qddot_Lagrangian}
\end{equation}
where the components of the matrices $J_{q,\dot{q}}$, $J_{\dot{q},\dot{q}}$ are given by
\begin{xalignat}{2}
  \left(J_{q,\dot{q}}\right)_{ij} &= \frac{\partial^2 \lag_{\text{NN}}}{\partial q_i \partial \dot{q}_j}, &
  \left(J_{\dot{q},\dot{q}}\right)_{ij} &= \frac{\partial^2 \lag_{\text{NN}}}{\partial \dot{q}_i \partial \dot{q}_j}.\label{eqn:J_matrices}
\end{xalignat}
The crucial point is now that the Lagrangian is constructed such that it is exactly invariant under the symmetry transformations $h^s$ which correspond to the conserved quantities. To guarantee that $\lag_{\text{NN}}$ satisfies Eq. \eqref{eqn:L_invariance}, we write $\lag_{\text{NN}} = D_{L} \circ D_{L-1} \circ D_1 \circ S$, where $D_\ell$ are standard, non-linear dense layers. The symmetry-enforcing first layer $S:\mathbb{R}^{2d}\rightarrow \mathbb{R}^{M}$ is a function that combines the inputs $q,\dot{q}$ into a set of $M$ invariants such that
\begin{equation}
  S_m(h^s(q),h^s_{*,q}(\dot{q})) = S_m(q,\dot{q})\qquad\text{for $m=1,2,\dots,M$}.\label{eqn:S_invariance}
\end{equation}
For example, if we consider a three-dimensional system which is invariant under the rotations in Eq. \eqref{eqn:rotations}, then the simplest three invariants that can be constructed from the position $q$ and the velocity $\dot{q}$ are $S_1(q,\dot{q})=q^2$, $S_2(q,\dot{q})=\dot{q}^2$ and $S_3(q,\dot{q})=q\cdot \dot{q}$; here and in the following we use the notation $u\cdot v := \sum_{j=1}^{d} u_jv_j$ to denote dot products of vectors $u,v\in\mathbb{R}^d$ and write $u^2:= u\cdot u$ for simplicity.

\begin{figure}
  \begin{center}
    \includegraphics[width=0.8\linewidth]{\figdir/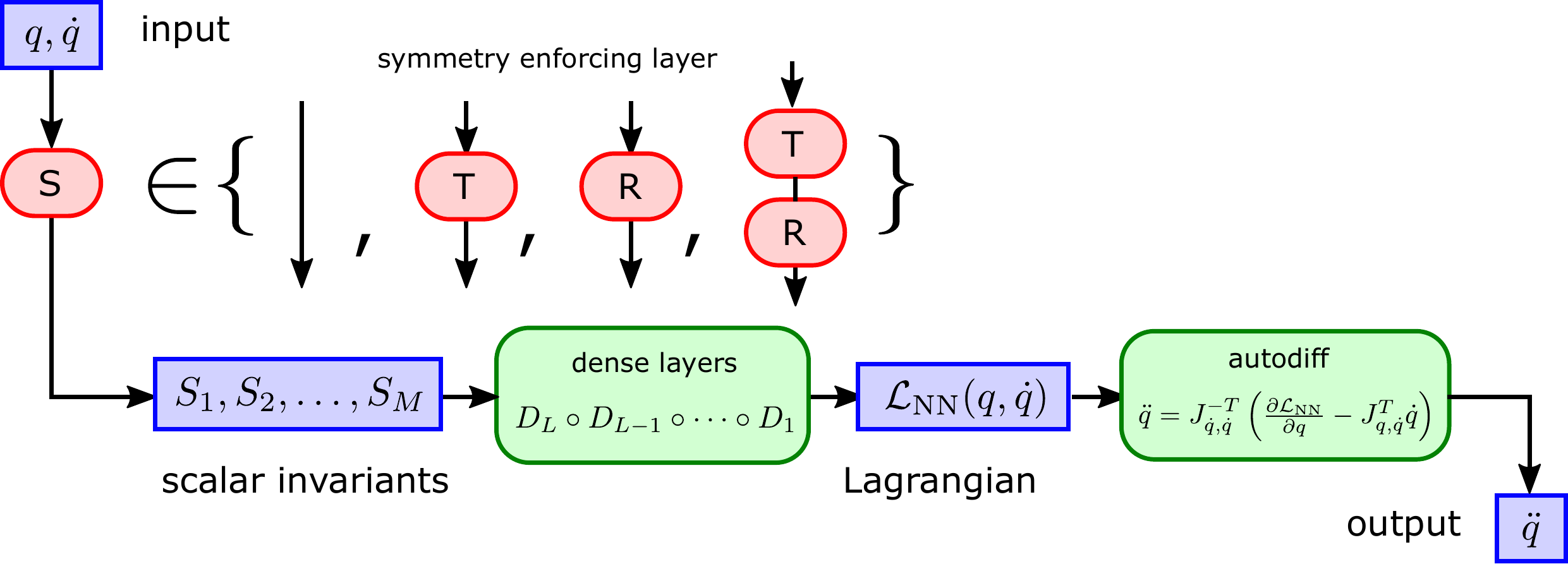}
    \caption{Structure of the function $\Phi$ which maps the inputs $q,\dot{q}$ to the predicted acceleration $\ddot{q}=\Phi(q,\dot{q})$. The symmetry-enforcing layers $T$, $R$ and $R\circ T$ for the specific systems studied in this work are discussed in the main text. A network without any symmetry constraints can be obtained by removing $S$, as represented by an undecorated downward arrow $\downarrow$.}
    \label{fig:nn_architecture}
  \end{center}
\end{figure}

The structure of the entire function $\Phi$ is shown in Fig. \ref{fig:nn_architecture}. Removing the symmetry-enforcing layer $S$ results in the standard architecture already introduced in \cite{Greydanus2019}. The explicit form of symmetry-enforcing layers for specific systems is discussed in Section \ref{sec:model_systems}.

The weights of the dense layers $D_\ell$ are learned by training the neural network on measurements of the physical system. As already pointed out above, the training data could consist of (noisy) observations of a real physical system. However, to demonstrate the principle and to not be restricted by a lack of training data, here we use synthetic, simulated data instead. More specifically, we compare the output of the neural network $\Phi_{\text{NN}}$ to (noisy) ground truth data. The ground truth predictions $\widehat{y}$ are generated by adding random noise to the true acceleration $\ddot{q} = \Phi_{\text{true}}(q,\dot{q})$ for some choice of $q$, $\dot{q}$; the function $\Phi_{\text{true}}$ is obtained by replacing $\lag_{\text{NN}}$ with the true Lagrangian $\lag_{\text{true}}$ in Eqs. \eqref{eqn:qddot_Lagrangian} and \eqref{eqn:J_matrices}. For each prediction $\widehat{y}$ the input $X$ of the neural network is obtained by also perturbing the corresponding $q,\dot{q}$ with random noise, thus simulating measurement errors for a real physical system.

During training we minimise the standard mean square error (MSE) loss function, which can be written as
\begin{equation}
  \text{Loss}\left(\{X^{(j)},\widehat{y}^{(j)}\}_{j=1}^{N}\right) = \frac{1}{N}\sum_{j=1}^{N}\left(\Phi(X^{(j)})-\widehat{y}^{(j)}\right)^2\label{eqn:MSE_loss}
\end{equation}
for a set of $N$ training samples $(X^{(1)},\widehat{y}^{(1)}),(X^{(2)},\widehat{y}^{(2)}),\dots,(X^{(N)},\widehat{y}^{(N)})$ with $X^{(j)}:=(q^{(j)}+\nu_1^{(j)},\dot{q}^{(j)}+\nu_2^{(j)})$ and $\widehat{y}^{(j)}=\ddot{q}^{(j)}+\nu_3^{(j)}$ where $\nu_i^{(j)}$ represents the random noise on the data.
\subsection{Model systems}\label{sec:model_systems}
We now discuss the three physical model systems, their symmetries and the corresponding conserved quantities in more detail.
\subsubsection{Motion of single particle in a gravitational potential}
\label{sec:kepler_problem}
First, we consider the dynamics of a particle of mass $m$ moving in a central force field with a potential energy that is inversely proportional to the distance from the origin. The true Lagrangian is
\begin{equation}
  \lag_{\text{true}}(q,\dot{q}) = \frac{m}{2}\dot{q}^2 + \frac{\alpha}{|q|}
  \qquad\text{with $q,\dot{q}\in\mathbb{R}^3$}
  \label{eqn:true_lagrangian_single_particle}
\end{equation}
for some positive constant $\alpha>0$. It is well known that for negative total energy $E=\frac{m}{2}\dot{q}^2-\alpha/|q|<0$ the trajectory is an ellipse where the origin of the coordinate system coincides with one focal point. Kepler studied this kind of motion in the solar system and derived three laws from observations, without knowing the analytical relationship between position, velocity and acceleration. The Lagrangian in Eq. \eqref{eqn:true_lagrangian_single_particle} is invariant under the rotations in Eq. \eqref{eqn:rotations}, which leads to the conservation of the vector-valued angular momentum $L$ in Eq. \eqref{eqn:angular_momentum_3d}. As a consequence, the entire trajectory lies in the plane spanned by the initial position $q(0)$ and initial velocity $\dot{q}(0)$; this plane is orthogonal to the vector $L$.

Our objective is to learn the neural network Lagrangian $\lag_{\text{NN}}(q,\dot{q})$ and hence the relationship between $q$, $\dot{q}$ and $\ddot{q}$ from data under the assumption that the dynamics is invariant under the action of the special orthogonal group $SO(3)$. To achieve this, recall that the simplest invariant scalars that can be constructed from position $q\in\mathbb{R}^3$ and velocity $\dot{q}\in\mathbb{R}^3$ are $q^2$, $\dot{q}^2$ and $q\cdot\dot{q}$. Hence, the symmetry-enforcing layer $S$ in Fig. \ref{fig:nn_architecture} should take as input $q$ and $\dot{q}$ and output these three invariant quantities; in the following we denote this layer as $R$ (for ``rotationally invariant'').

Adopting the Einstein sum convention of summing over pairs of identical upper and lower indices, the three quantities that are conserved for the dynamics generated with neural network Lagrangian $\lag_{\text{NN}}$ can be written as
\begin{equation}
  L^{(\rho)}_{\text{NN}} = \varepsilon^{\rho\sigma\tau} q_\sigma\frac{\partial \lag_{\text{NN}}}{\partial \dot{q}_\tau}
  \label{eqn:angular_momentum_NN}
\end{equation}
where $\varepsilon$ is the  three-dimensional Levi Civita symbol defined in  Eq. \eqref{eqn:levi_civita_definition}. We stress that by construction the angular momenta in Eq. \eqref{eqn:angular_momentum_NN} are \textit{exactly} conserved, independent of the weights of the neural network.
The quantity $L^{(\rho)}_{\text{NN}}$ can be likened to the shadow Hamiltonian, which is conserved for symplectic integrators of time-independent Hamiltonian systems. For the true Lagrangian in Eq. \eqref{eqn:true_lagrangian_single_particle} the ``neural network'' angular momentum $L^{(\rho)}_{\text{NN}}$ reduces to the three components of the ``true'' angular momentum vector in Eq. \eqref{eqn:angular_momentum_3d}, namely
\begin{equation}
  L^{(\rho)}_{\text{true}} = m \varepsilon^{\rho\sigma\tau} q_\sigma\dot{q}_\tau.
  \label{eqn:angular_momentum_true}
\end{equation}
For future reference we combine the quantities defined in Eqs. \eqref{eqn:angular_momentum_NN} and \eqref{eqn:angular_momentum_true} into two vectors:
\begin{xalignat}{2}
  L_{\text{NN}} &= \left(L_{\text{NN}}^{(1)},L_{\text{NN}}^{(2)},L_{\text{NN}}^{(3)}\right),&
  L_{\text{true}} &= \left(L_{\text{true}}^{(1)},L_{\text{true}}^{(2)},L_{\text{true}}^{(3)}\right).\label{eqn:angular_momentum_kepler_vector}
\end{xalignat}
\subsubsection{Motion of a massive relativistic particle in the Schwarzschild metric}\label{sec:schwarzschild_problem}
Next, we consider the motion of a massive particle in the Schwarzschild metric \cite{Schwarzschild1916,Droste1917}, which is the relativistic equivalent to the Kepler problem in Section \ref{sec:kepler_problem}. With the space-time dependent metric tensor $g(q)$ the Lagrangian can be written as\footnote{The geodesic line element is $ds=\sqrt{g_{\mu\nu}dq^\mu dq^\nu}$, which would imply that the Lagrangian is $ds/dt=\sqrt{g_{\mu\nu}(q)\dot{q}^\mu \dot{q}^\nu}$. However, since the square root is a monotonous function, the Lagrangian in Eq. \eqref{eqn:schwarzschild_lagrangian} has the same stationary points and thus generates the same dynamics.}
\begin{equation}
  \lag = g_{\mu\nu}(q)\dot{q}^\mu \dot{q}^\nu \qquad \text{for $\mu,\nu=1,2,3,4$},\label{eqn:schwarzschild_lagrangian}
\end{equation}
where $q^\mu$ are the components of the four-dimensional contravariant time-position vector $q=(x,\tau)=(x_1,x_2,x_3,\tau)$. $\dot{q}^\mu = dq^\mu/dt$ denotes the derivative with respect to the eigen-time $t$ experienced by a moving observer and $\tau$ is the time measured by a static observer far away from the origin, where $g \rightarrow \operatorname{diag}{(+1,+1,+1,-1)}$ tends to the constant metric of flat space-time.

Expressing the spatial coordinate $x$ in spherical coordinates, the true Lagrangian in Eq. \eqref{eqn:schwarzschild_lagrangian} can be written as
\begin{equation}
  \lag_{\text{true}} = -\left(1-\frac{r_s}{r}\right)\dot{\tau}^2 + \left(1-\frac{r_s}{r}\right)^{-1}\dot{r}^2 + r^2\left(\dot{\theta}^2 + \sin^2(\theta)\dot{\varphi}^2\right) \label{eqn:schwarzschild_lagrangian_II}
\end{equation}
where $r_s$ is the Schwarzschild radius and $r=|x|$ is the distance from the origin. It is easy to see\footnote{Observe that $r^2\left( \dot{\theta}^2 + \sin^2(\theta)\dot{\varphi}^2\right) = \dot{x}^2 - (x\cdot\dot{x})^2/x^2$.} that Eq. \eqref{eqn:schwarzschild_lagrangian_II} is invariant under rotations in three dimensional space and can be expressed entirely in terms of the three scalars $x^2$, $\dot{x}^2$, $x\cdot\dot{x}$ and the time derivative $\dot{\tau}$. As a consequence, the dynamics conserves the (specific) three-dimensional angular momentum which turns out to be
\begin{equation}
  L = x\times \dot{x}.
\end{equation}
Note that this expression differs from the one in Eq. \eqref{eqn:angular_momentum_3d} only through scaling by the (constant) mass $m$. In analogy to Eq. \eqref{eqn:angular_momentum_true} we can also write down the three components of the angular momentum vector as
\begin{equation}
  L^{(\rho)}_{\text{true}} = \varepsilon^{\rho\sigma\tau} x_\sigma\dot{x}_\tau.
  \label{eqn:relativistic_angular_momentum_true}
\end{equation}
Again, we want to learn the Lagrangian $\lag_{\text{NN}}(q,\dot{q})$ which is represented by a neural network. We drop the dependency on $\tau$, since we found that including it leads to instabilities during training. The reason for this is that for the true solution $\tau$ grows (approximately linearly) with time. While the neural network will learn that the coefficient that multiplies $\tau$ in the Lagrangian is small, it will always remain non-zero, thus leading to a large, unphysical contribution for $\tau\gg 1$. Note that dropping $\tau$ from the inputs is consistent with the true Lagrangian in Eq. \eqref{eqn:schwarzschild_lagrangian_II}, which only depends on $\dot{\tau}$ but not the  $\tau$ itself. The symmetry enforcing layer takes as input the two vectors $q,\dot{q}\in\mathbb{R}^4$ and it will return the four scalars $\dot{\tau}$, $x^2$, $\dot{x}^2$ and $x\cdot\dot{x}$.

With this symmetry enforcing layer the dynamics generated by the neural network Lagrangian conserves the three quantities
\begin{equation}
  L^{(\rho)}_{\text{NN}} = \varepsilon^{\rho\sigma\tau}x_\sigma\frac{\partial \lag_{\text{NN}}}{\partial \dot{x}_\tau}
  \qquad\text{for $\rho=1,2,3$}\label{eqn:relativistic_angular_momentum_NN}
\end{equation}
exactly.
\subsubsection{Two interacting particles in $D$ dimensions}\label{sec:two_particle_problem}
Finally, we consider a system of two non-relativistic particles with masses $m_1$ and $m_2$ that move in $D$-dimensional space and interact via a potential that is invariant under translations and rotations. Setting $d=2D$, the $d$-dimensional state vector $q=(x^{(1)},x^{(2)})\in\mathbb{R}^d$ contains the two particle positions $x_1,x_2\in\mathbb{R}^D$ with corresponding velocities $\dot{x}^{(2)},\dot{x}^{(2)}\in\mathbb{R}^D$. The true Lagrangian is
\begin{equation}
  \lag_{\text{true}}(q,\dot{q}) = \lag_{\text{true}}(x^{(1)},x^{(2)},\dot{x}^{(1)},\dot{x}^{(2)})
  = \frac{m_1}{2}(\dot{x}^{(1)})^2 + \frac{m_2}{2}(\dot{x}^{(1)})^2 - V(||x^{(1)}-x^{(2)}||_2)\label{eqn:true_lagrangian_two_particle}
\end{equation}
with the interaction potential $V(r)$ given by the double-well function
\begin{equation}
  V(r) = \frac{\mu}{2}r^2 - \frac{\kappa}{4}r^4\label{eqn:double_well}
\end{equation}
for some positive constants $\mu,\kappa>0$. As usual, $||\cdot||_2$ denotes the Euclidean two-norm
\begin{equation}
  ||z||_2 := \left(\sum_{j=1}^{D}z_j^2\right)^{1/2}\qquad\text{for $\in\mathbb{R}^D$}.
\end{equation}
Again, the goal is to learn the neural network Lagrangian
\begin{equation}
  \lag_{\text{NN}}(q,\dot{q}) = \lag_{\text{NN}}(x^{(1)},x^{(2)},\dot{x}^{(1)},\dot{x}^{(2)})
  \label{eqn:nn_lagrangian_two_particle}
\end{equation}
which approximates the dynamics of the system, while taking into a account the symmetries of the problem. Although $\lag_{\text{NN}}$ in Eq. \eqref{eqn:nn_lagrangian_two_particle} depends on $4D$ unknowns, it can be restricted considerably by assuming that the system is invariant under translations and/or rotations. For this, we consider the following two continuous transformations of the entire system:
\begin{description}
  \item[Translations] $h^s_T$ with a constant vector $\Delta\in\mathbb{R}^D$
        \begin{equation}
          h^s_T: q = (x^{(1)}, x^{(2)})\mapsto (x^{(1)}+s\Delta,x^{(2)}+s\Delta)
        \end{equation}
  \item[Rotations] $h^s_R$ with a constant rotation matrix $R=\exp[s\Gamma]\in SO(D)$
        \begin{equation}
          h^s_R: q = (x^{(1)}, x^{(2)})\mapsto (R(s)x^{(1)},R(s)x^{(2)}),
        \end{equation}
        where the antisymmetric $D\times D$ matrix $\Gamma$ belongs to the fundamental representation of the Lie algebra $\mathfrak{so}(D)$ of the special orthogonal group $SO(D)$.
\end{description}
Assuming that the physics of the system is invariant under translations and/or rotations dramatically restricts the possible form of the Lagrangian $\lag_{\text{NN}}$. For example, for a rotationally invariant (but not necessarily translationally invariant) system the Lagrangian can only depend on dot-products of pairs of the four dynamical variables $x^{(1)},x^{(2)},\dot{x}^{(1)},\dot{x}^{(2)}\in\mathbb{R}^D$ and contractions of these variables with the antisymmetric Levi-Civita symbol $\varepsilon$ in $D$ dimensions. More generally, for a set of $n$ vectors $A=\{a^{(1)},a^{(2)},\dots,a^{(n)}\}$ with $a^{(j)}\in\mathbb{R}^D$ we denote by $\mathscr{R}(A)$ the set of all rotationally invariant scalars; a detailed discussion of the construction of this set can be found in \ref{sec:rotationally_invariant_combinations} (see also \cite[Section 3, Lemma 2]{Villar2021}). For example we have for the three vectors $u,v,w\in\mathbb{R}^3$:
\begin{equation}
  \mathscr{R}(\{u,v,w\}) = \{u^2,v^2,w^2,u\cdot v,u\cdot w,v\cdot w,u\cdot(v\times w)\}
\end{equation}
since there are six different scalar products and
\begin{equation}\varepsilon^{ijk}u_iv_jw_k = u\cdot (v\times w)
\end{equation}
is the only non-vanishing contraction with the three-dimensional Levi-Civita symbol in this case.
With this notation, the set of scalar invariants that are output by the symmetry-enforcing layer $S$ can be written down as in Tab. \ref{tab:lagrangian_restrictions}.
\begin{table}
  \begin{center}
    \begin{tabular}{|cc|c|c|}
      \hline
                      &     & \multicolumn{2}{c|}{rotationally invariant?}                                                                                                     \\
                      &     & no                                                               & yes                                                                           \\
      \hline
      translationally & no  & $\{x^{(1)},x^{(2)},\dot{x}^{(1)},\dot{x}^{(2)}\}$                & $R(q,\dot{q}) = \mathscr{R}(\{x^{(1)},x^{(2)},\dot{x}^{(1)},\dot{x}^{(2)}\})$ \\
      \cline{2-4}
      invariant?      & yes & $T(q,\dot{q}) = \{x^{(1)}-x^{(2)},\dot{x}^{(1)},\dot{x}^{(2)}\}$ &
      $(R\circ T)(q,\dot{q}) = \mathscr{R}(\{x^{(1)}-x^{(2)},\dot{x}^{(1)},\dot{x}^{(2)}\})$                                                                                   \\
      \hline
    \end{tabular}
    \caption{Output of the symmetry-enforcing layer $S$.}
    \label{tab:lagrangian_restrictions}
  \end{center}
\end{table}
\begin{figure}
  \begin{center}
    \includegraphics[width=0.9\linewidth]{\figdir/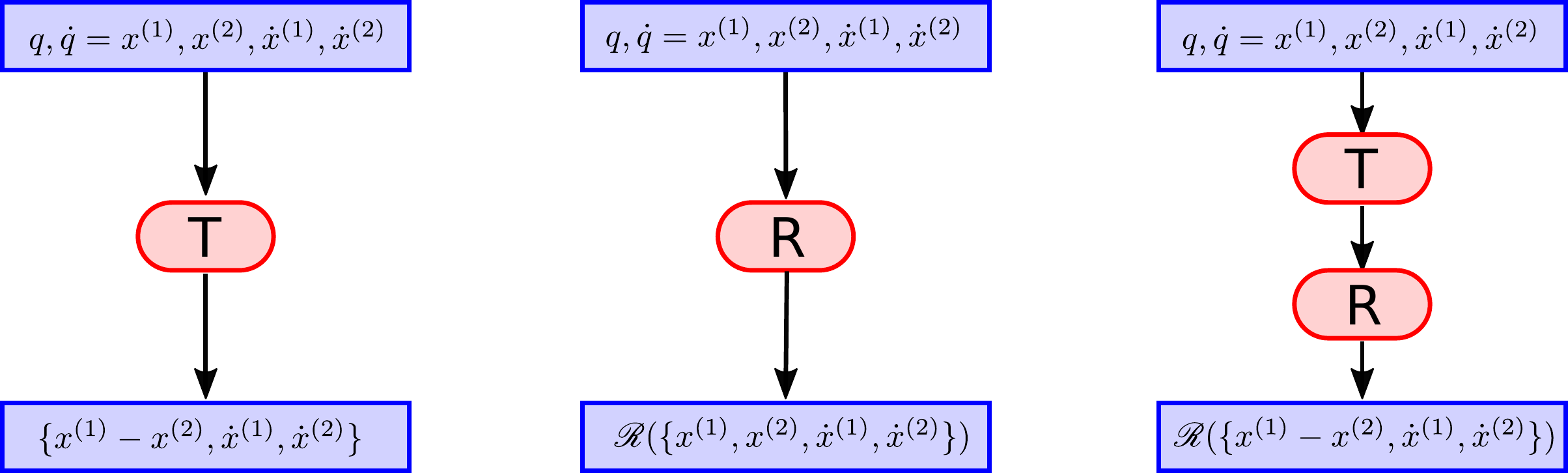}
    \caption{Pictorial representation of the symmetry enforcing layers.}
    \label{fig:symmetry_enforcing_layers}
  \end{center}
\end{figure}
While the most general Lagrangian is a function of $4D$ unknowns, translational invariance reduces it to a function which only depends on $3D$ variables. A Lagrangian which is invariant under rotations (but not necessarily translations) depends on $10+{4\choose D}$ unknowns; assuming both translational and rotational invariance reduces this further to only $6+{3\choose D}$ variables.

To derive the conserved quantities according to Eq. \eqref{eqn:noether_conserved_I}, observe that the $D$ generators of the translation group are the vectors $\Delta^{(\alpha)}\in\mathbb{R}^D$ with components
\begin{equation}
  \Delta^{(\alpha)}_j = \delta_{\alpha j}=\begin{cases}
    1 & \text{if $\alpha=j$} \\
    0 & \text{otherwise}
  \end{cases}
  \qquad\text{for $1\le\alpha\le D$}
\end{equation}
where $\delta_{jk}$ is the Kronecker-$\delta$. The $\frac{1}{2}D(D-1)$ generators $\Gamma^{(\rho,\sigma)}$ of the rotation group are the antisymmetric $D\times D$ matrices with entries
\begin{equation}
  \Gamma^{(\rho,\sigma)}_{jk} = \delta_{\rho j}\delta_{\sigma k} - \delta_{\rho k}\delta_{\sigma j}\qquad\text{for $1\le\rho<\sigma\le D$}.
\end{equation}
Using Noether's Theorem, this leads to the following $\frac{1}{2}D(D+1)$ conserved quantities $M^{(\alpha)}_{\text{NN}}$ and $J^{(\rho,\sigma)}_{\text{NN}}$:
\begin{equation}
  \begin{aligned}
    M^{(\alpha)}_{\text{NN}}      & = \frac{\partial \lag_{\text{NN}}}{\partial \dot{x}^{(1)}_\alpha} + \frac{\partial \lag_{\text{NN}}}{\partial \dot{x}^{(2)}_\alpha}\qquad\text{for $1\le\alpha\le D$},                                                                                                                                                                                            \\
    L^{(\rho,\sigma)}_{\text{NN}} & = \frac{\partial \lag_{\text{NN}}}{\partial \dot{x}^{(1)}_\rho} x_\sigma^{(1)}-\frac{\partial \lag_{\text{NN}}}{\partial \dot{x}^{(1)}_\sigma} x_\rho^{(1)} + \frac{\partial \lag_{\text{NN}}}{\partial \dot{x}^{(2)}_\rho} x_\sigma^{(2)}-\frac{\partial \lag_{\text{NN}}}{\partial \dot{x}^{(2)}_\sigma} x_\rho^{(2)}\qquad\text{for $1\le\rho < \sigma\le D$}.
  \end{aligned}
  \label{eqn:two_particle_NN_momenta}
\end{equation}
For the true Lagrangian in Eq. \eqref{eqn:true_lagrangian_two_particle} this simplifies to the usual linear momenta $M^{(\alpha)}_{\text{true}}$ and (generalised) angular momenta $L^{(\rho,\sigma)}_{\text{true}}$ with
\begin{equation}
  \begin{aligned}
    M^{(\alpha)}_{\text{true}}      & = m_1\dot{x}^{(1)}_\alpha + m_2\dot{x}^{(2)}_\alpha\qquad\text{for $1\le\alpha\le D$},                                                                                                                                \\
    L^{(\rho,\sigma)}_{\text{true}} & = m_1\left(\dot{x}^{(1)}_\rho x_\sigma^{(1)}-\dot{x}^{(1)}_\sigma x_\rho^{(1)}\right) + m_2\left(\dot{x}^{(2)}_\rho x_\sigma^{(2)}-\dot{x}^{(2)}_\sigma x_\rho^{(2)}\right)\qquad\text{for $1\le\rho < \sigma\le D$}.
  \end{aligned}\label{eqn:momentum_true_twoparticle}
\end{equation}
For future reference we collect the conserved quantities in the following vectors
\begin{equation}
  \begin{aligned}
    M_{\text{NN}}   & = \left(M_{\text{NN}}^{(1)},M_{\text{NN}}^{(2)},\dots,M_{\text{NN}}^{(D)}\right)^T\in\mathbb{R}^D,                                 \\
    L_{\text{NN}}   & = \left(L_{\text{NN}}^{(1,2)},L_{\text{NN}}^{(1,3)},\dots,M_{\text{NN}}^{(D-1,D)}\right)^T\in\mathbb{R}^{\frac{1}{2}D(D-1)},       \\
    M_{\text{true}} & = \left(M_{\text{true}}^{(1)},M_{\text{true}}^{(2)},\dots,M_{\text{true}}^{(D)}\right)^T\in\mathbb{R}^D,                           \\
    L_{\text{true}} & = \left(L_{\text{true}}^{(1,2)},L_{\text{true}}^{(1,3)},\dots,M_{\text{true}}^{(D-1,D)}\right)^T\in\mathbb{R}^{\frac{1}{2}D(D-1)}.
  \end{aligned}
  \label{eqn:two_particle_momentum_vectors}
\end{equation}
\section{Results}\label{sec:results}
\subsection{Implementation and hyperparameters}\label{sec:implementation}
The code that was used to obtain the numerical results reported in this section has been implemented in tensorflow  \cite{Tensorflow2015} and is freely available at the following URL, which also contains instructions on how to install and run the code:
\begin{center}
  \url{https://github.com/eikehmueller/mlconservation_code}
\end{center}
The results reported in this paper were generated with the release published at \cite{code_release}. In all cases the Lagrangian is represented by two hidden dense layers $D_1$, $D_2$ with $n_{\text{h}}=128$ output units each and a final dense layer $D_3$ with a single output unit. This unit does not have a bias term since adding a constant to the Lagrangian does not change the equations of motion. A softplus activation function is used in all cases. The layer weights are initialised with a random normal distribution similarly to \cite{Greydanus2019}: the standard deviation of the normal distribution is set to $2/\sqrt{n_{\text{h}}}$ for the first hidden layer $D_1$, $1/\sqrt{n_{\text{h}}}$ for the second hidden layer $D_2$ and $\sqrt{n_{\text{h}}}$ for the output layer $D_3$; the biases of all layers are initialised to zero. Single precision arithmetic is used in all numerical experiments.

For all three systems discussed in Section \ref{sec:model_systems} the networks are trained over 2500 epochs with 100 steps per epoch and a batch size of $B=128$ using the Adam optimiser. The training schedule is a cosine decay, starting with a learning rate of $10^{-3}$ which is reduced to $10^{-5}$ at the end of the training cycle. We find that this reduces the MSE to $10^{-5}-10^{-6}$ for the single-particle problem described in Section \ref{sec:kepler_problem}, to $10^{-6}$ for the relativistic particle in Section \ref{sec:schwarzschild_problem} and to approximately $10^{-4}$ for the two-particle system in Section \ref{sec:two_particle_problem}. This is consistent with the chosen random noise on the training data, which is $\sigma=10^{-3}$ in all cases and therefore limits the minimal achievable MSE to the order of $\sigma^2\sim 10^{-6}$. The full loss histories can be found in \ref{sec:loss_histories}.
\subsubsection{Motion of a single particle in a gravitational potential}\label{sec:results_single_particle}
Both the mass of the particle and the strength of the gravitational potential in Eq. \eqref{eqn:true_lagrangian_single_particle} are set to $m=\alpha=1$. To train the model, we generate pairs of inputs $X^{(j)} = (q_{\text{exact}}(\varphi^{(j)})+\sigma \xi_1^{(j)},\dot{q}_{\text{exact}}(\varphi^{(j)})+\sigma\xi_2^{(j)})\in\mathbb{R}^6$ and ground truth $\widehat{y}^{(j)} = \ddot{q}_{\text{exact}}(\varphi^{(j)})+\sigma\xi_3^{(j)}\in\mathbb{R}^3$ where $q_{\text{true}}(\varphi)$ is the exact solution as a function of the angle azimuthal angle $\varphi$. The size of the noise is characterised by $\sigma=10^{-3}$ as discussed above. The angles $\varphi^{(j)}\sim\text{Uniform}(-\pi,+\pi)$ are uniformly distributed while $\xi^{(j)}_{1,2,3}\sim\mathcal{N}(0,1)$ are drawn from a normal distribution with mean zero and unit variance. We assume that the motion takes place entirely in the $x-y$ plane and the true trajectory is an ellipse with eccentricity $\varepsilon_{\text{ecc}}=0.8$. The initial conditions are chosen such that the vertical component of the angular momentum is $L_z=1$.

\begin{figure}
  \begin{center}
    \includegraphics[width=\linewidth]{\figdir/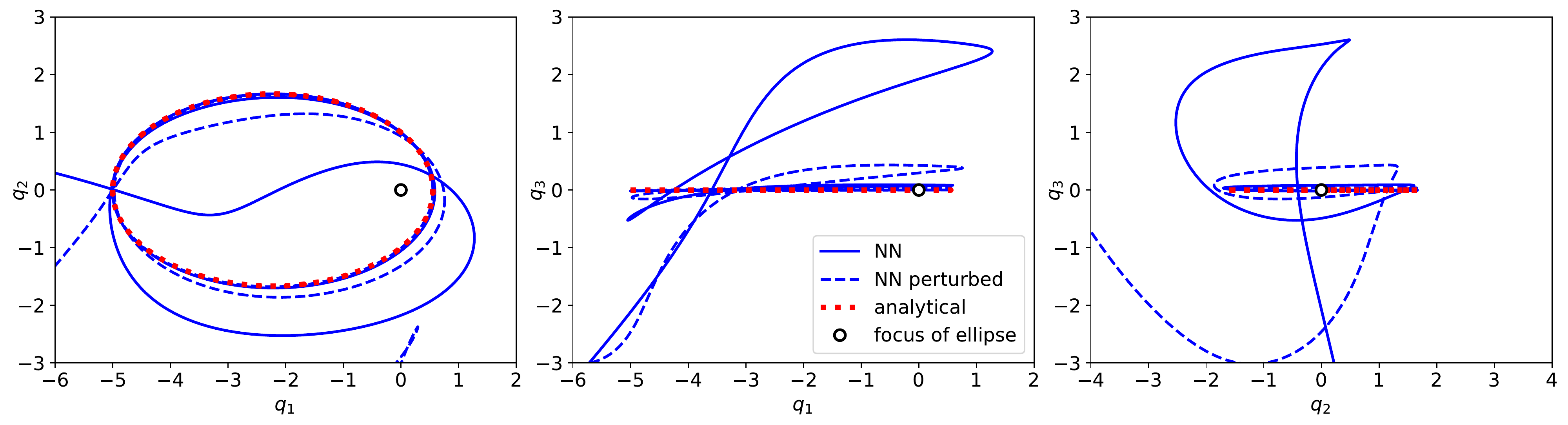}\\
    \includegraphics[width=\linewidth]{\figdir/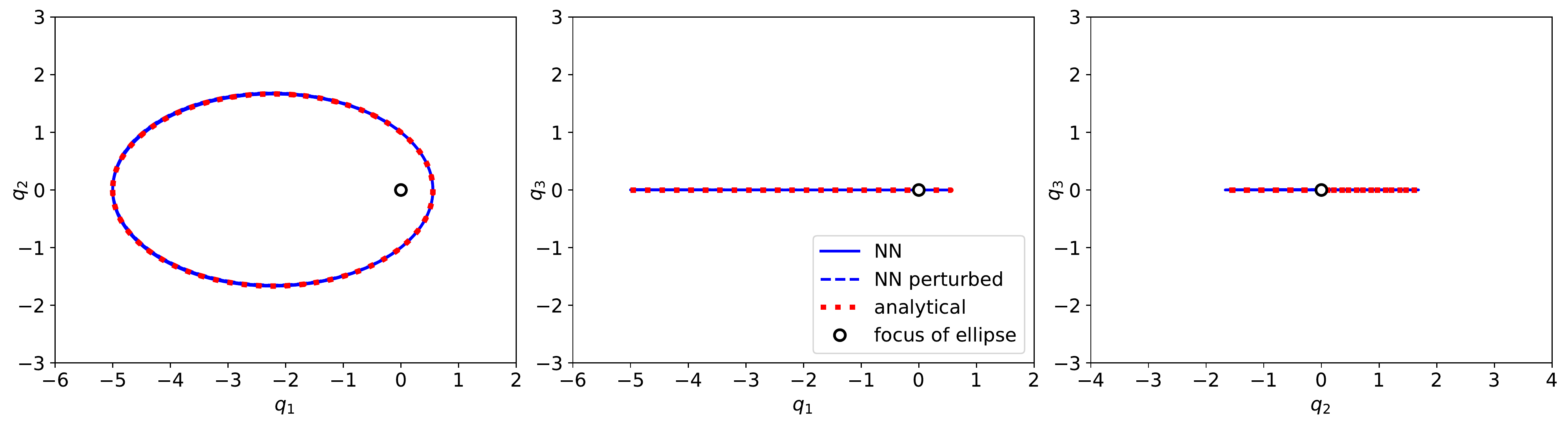}
    \caption{Trajectories for the motion of a particle in a gravitational potential described in Section \ref{sec:kepler_problem} without (top) and with (bottom) constraints on the neural network to enforce rotational invariance of the Lagrangian. In each case, a trajectory that is obtained by perturbing the initial conditions by $\sim10^{-3}$ is also shown as a dashed curve. The ellipse that represents the true solution is marked with red dots.}
    \label{fig:trajectories_kepler}
  \end{center}
\end{figure}

\begin{figure}
  \begin{center}
    \includegraphics[width=\linewidth]{\figdir/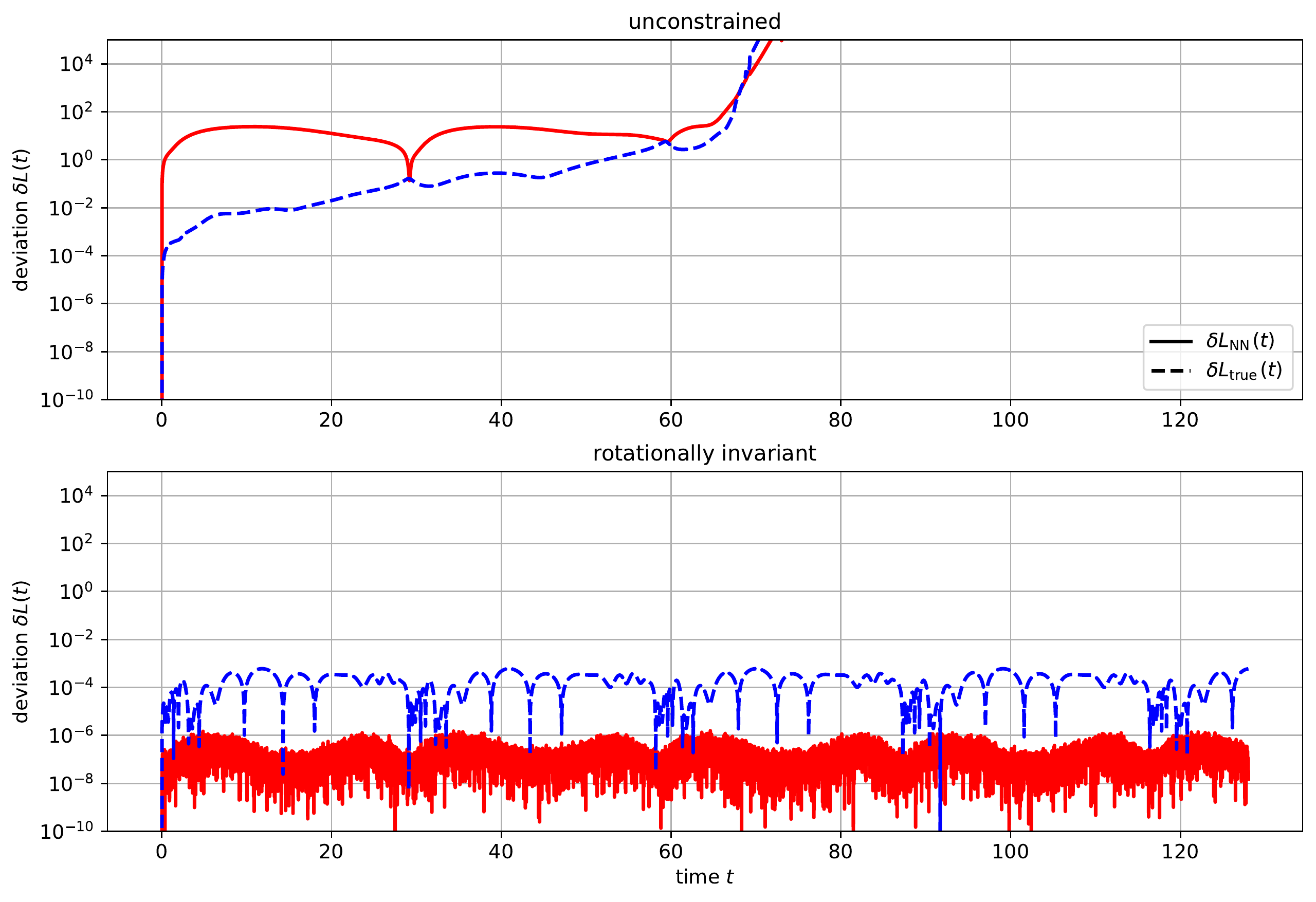}
    \caption{Conservation of angular momentum for     the motion of a particle in a gravitational potential defined in Section \ref{sec:kepler_problem} without (top) and with (bottom) constraints on the neural network to enforce rotational invariance of the Lagrangian. In each case the deviation $\delta L_{\text{NN}}(t)$ and $\delta L_{\text{true}}(t)$ of the absolute value of the total three dimensional angular momentum vectors are shown as a function of time. Both the ``neural network'' angular momentum $L_{\text{NN}}(t)$ defined in Eq. \eqref{eqn:angular_momentum_NN} and the true angular momentum $L_{\text{true}}(t)$ defined in Eq. \eqref{eqn:angular_momentum_true} are considered.}
    \label{fig:conservation_kepler}
  \end{center}
\end{figure}

\begin{figure}
  \begin{center}
    \includegraphics[width=\linewidth]{\figdir/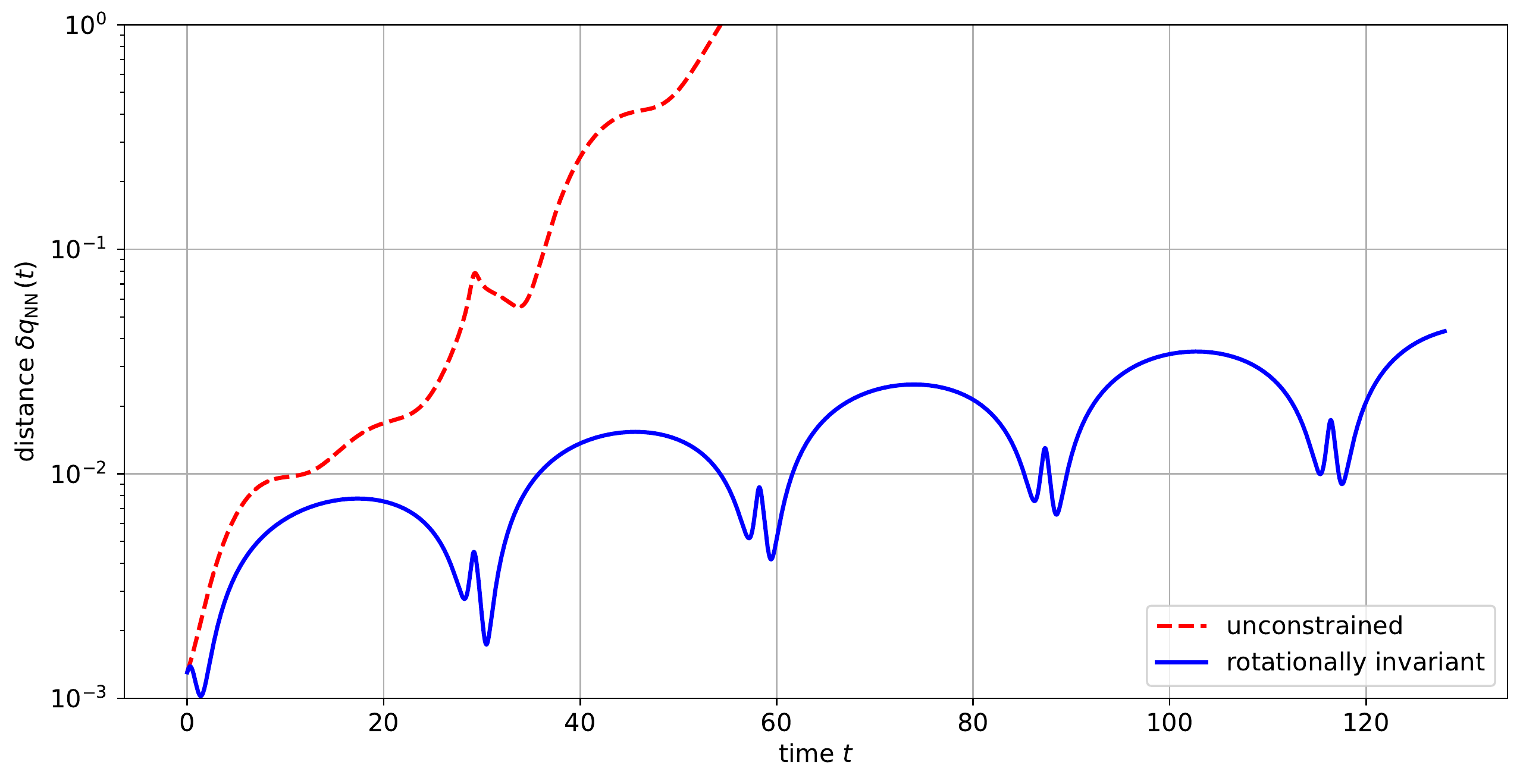}
    \caption{Evolution of the distance $\delta q_{\text{NN}}(t)= ||q_{\text{NN}}(t)-q_{\text{NN}}^{\text{(perturbed)}}(t)||_2$ between two trajectories obtained with slightly different initial conditions for the Kepler problem. The initial conditions that differ by $\delta q_{\text{NN}}(0)\sim 10^{-3}$.}
    \label{fig:pertubation_kepler}
  \end{center}
\end{figure}

Fig. \ref{fig:trajectories_kepler} shows the trajectories predicted with the trained neural network Lagrangian. The equations of motion are integrated up to the final time $T=128$ with a fourth order Runge Kutta (RK4) method and a timestep size of $\Delta t=10^{-2}$. The figure also shows a trajectory for which both initial position and initial velocity are perturbed by normal random noise with a standard deviation of $10^{-3}$. As can be seen from the top figure, the trajectory obtained with the neural network Lagrangian deviates significantly from the true solution (the red ellipse) for the unconstrained Lagrangian and the two trajectories with perturbed initial conditions diverge. In fact, the neural network trajectories become unstable. Radically different behaviour is observed for the neural network with built-in rotational invariance (bottom plot in Fig. \ref{fig:trajectories_kepler}). Visually the trajectories obtained with the neural network Lagrangian can not be distinguished from the true solution. In particular, the trajectories appear to be confined to the $x-y$ plane, which implies that the two corresponding components of the angular momentum indeed remain close to zero.
To investigate the conservation of angular momentum further, we compute the time evolution of the deviation of the angular momentum from its original value. For this we consider the deviation
\begin{xalignat}{2}
  \delta L_{\text{NN}}(t) & = \frac{\left|\left|L_{\text{NN}}(t)-L_{\text{NN}}(0)\right|\right|_2}{\left|\left|L_{\text{NN}}(0)\right|\right|_2}, &
  \delta L_{\text{true}}(t) & =\frac{\left|\left|L_{\text{true}}(t)-L_{\text{true}}(0)\right|\right|_2}{\left|\left|L_{\text{true}}(0)\right|\right|_2}.
  \label{eqn:deltaL_definition}
\end{xalignat}
of the angular momentum vectors defined by Eqs. \eqref{eqn:angular_momentum_NN}, \eqref{eqn:angular_momentum_true} and \eqref{eqn:angular_momentum_kepler_vector} from their initial values at time $t=0$. Exact conservation of angular momentum would correspond to $\delta L_{\text{true}}(t) = 0$ for all times $t$. As Fig. \ref{fig:conservation_kepler} (top) shows, for the unconstrained network $\delta L_{\text{true}}(t)$ is small initially, since the network has learned some degree of angular momentum conservation from the data, but then increases to around $1$ by time $t\approx 50$ and diverges shortly after that. As expected, for the rotationally invariant network shown in the bottom figure $\delta L_{\text{NN}}(t)$ remains zero within single precision rounding errors. More importantly, for the rotationally invariant neural network the deviation $\delta L_{\text{true}}(t)$ of the true angular momentum never exceeds $10^{-3}$. Hence, even though mathematically only the conservation of $L_{\text{NN}}$ can be guaranteed, numerically this appears to also help with the conservation of $L_{\text{true}}$. This behaviour is similar to the approximate conservation of the true energy observed for symplectic integrators of time-independent Hamiltonian systems, which is related to the exact conservation of the shadow Hamiltonian.
Finally, in Fig. \ref{fig:pertubation_kepler} we show the distance $\delta q_{\text{NN}}(t) = ||q_{\text{NN}}(t)-q_{\text{NN}}^{(\text{perturbed})}(t)||_2$ between the position vectors of the two neural network trajectories that only differ by a $\sim10^{-3}$ perturbation of the initial condition. While -- as expected -- the trajectories diverge over time, for the  rotationally invariant neural network they stay much closer together and for the considered time interval their distance stays below $0.1$. Intuitively, the reason for this is that rotational invariance limits the allowed trajectories to a smaller sub-manifold, so there is less ``room'' for nearby trajectories to diverge.
\subsubsection{Motion of a massive relativistic particle in the Schwarzschild metric}
For the motion of a massive relativistic particle in the Schwarzschild metric as defined in Section \ref{sec:schwarzschild_problem} we set the Schwarzschild radius to $r_s=0.1$. In this case, the ``true'' trajectory $q_{\text{exact}}(t)$ is obtained by integrating the exact equations of motion, which can be obtained from the true Lagrangian in Eq. \eqref{eqn:schwarzschild_lagrangian_II}, with a RK4 integrator with a timestep size of $\Delta t=10^{-2}$. As above, normally distributed random noise with a standard deviation of $\sigma=10^{-3}$ is added to obtain training samples $(X^{(j)},\widehat{y}^{(j)})$ with $X^{(j)} = (q_{\text{exact}}(t^{(j)})+\sigma \xi_1^{(j)},\dot{q}_{\text{exact}}(t^{(j)})+\sigma\xi_2^{(j)})\in\mathbb{R}^8$ and ground truth $\widehat{y}^{(j)} = \ddot{q}_{\text{exact}}(t^{(j)})+\sigma\xi_3^{(j)}\in\mathbb{R}^4$, where $\xi_{1}^{(j)}, \xi_{2}^{(j)}, \xi_{3}^{(j)}\sim \mathcal{N}(0,1)$ and $t^{(j)}$ are sample times along the true trajectory.

\begin{figure}
  \begin{center}
    \includegraphics[width=\linewidth]{\figdir/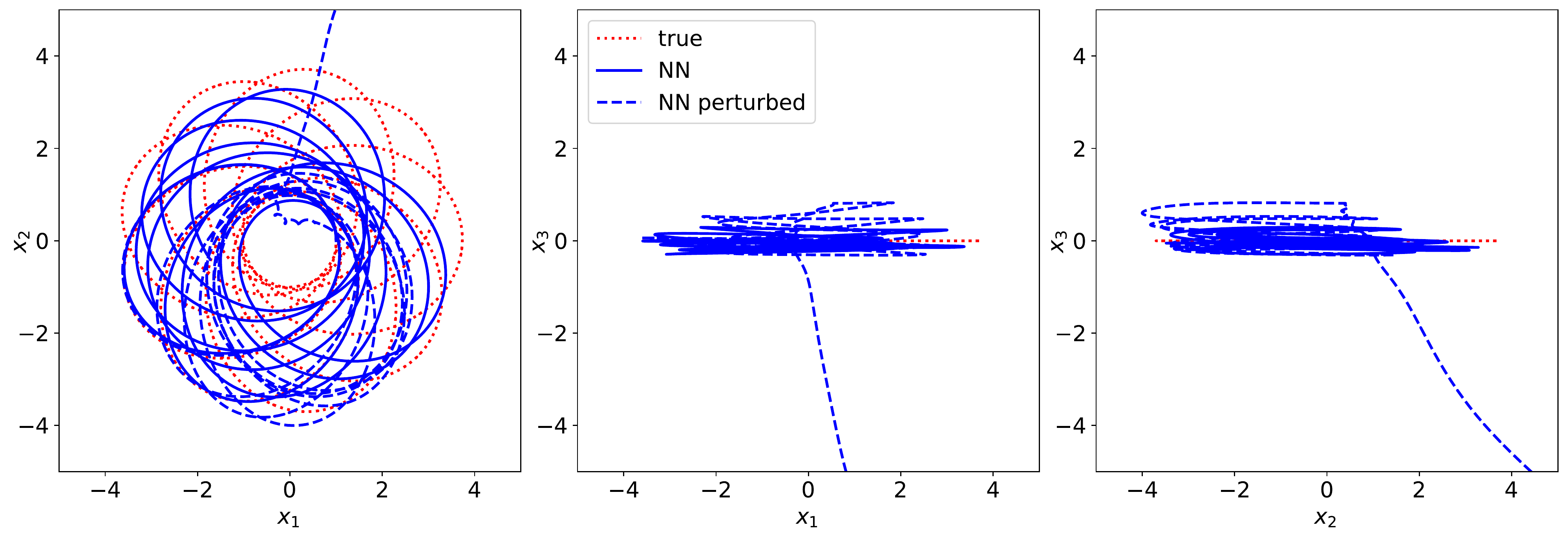}\\
    \includegraphics[width=\linewidth]{\figdir/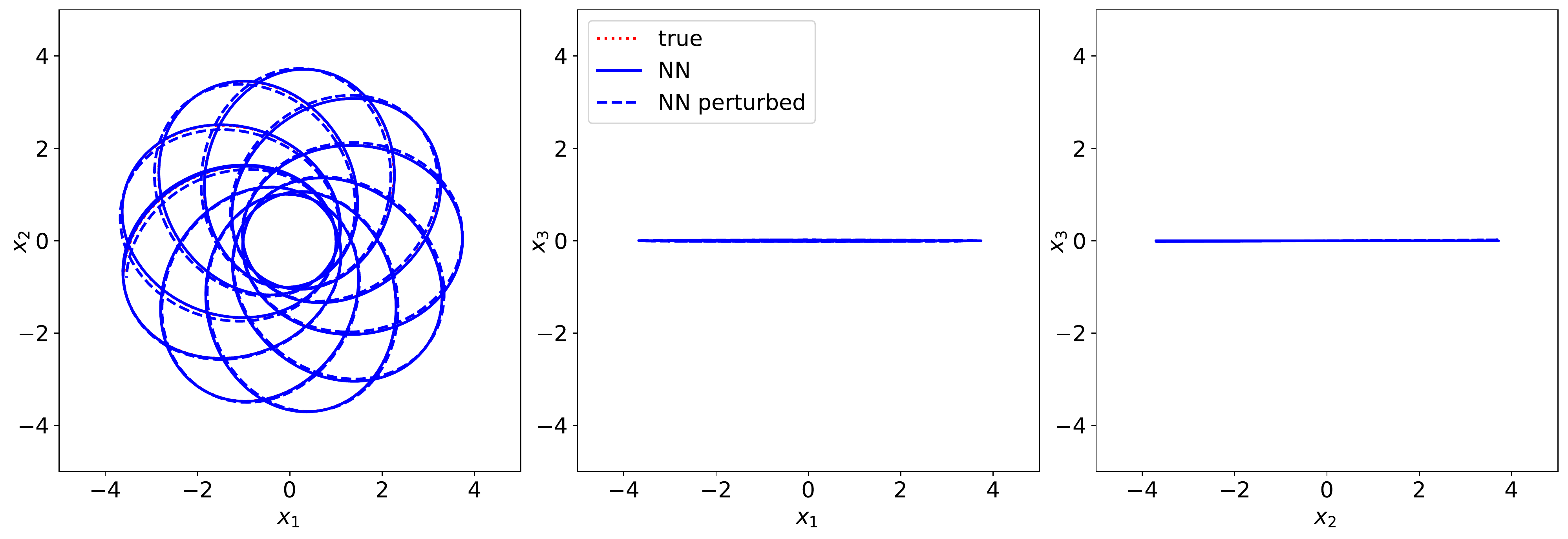}
    \caption{Trajectories for the motion of a relativistic massive particle in the Schwarzschild metric as defined in Section \ref{sec:schwarzschild_problem} without (top) and with (bottom) constraints on the neural network to enforce rotational invariance of the spatial part of the Lagrangian. In each case, a trajectory that is obtained by perturbing the initial conditions by $\sim10^{-3}$ is also shown as a dashed curve. The true solution is shown as a dotted red curve.}
    \label{fig:trajectories_schwarzschild}
  \end{center}
\end{figure}

\begin{figure}
  \begin{center}
    \includegraphics[width=\linewidth]{\figdir/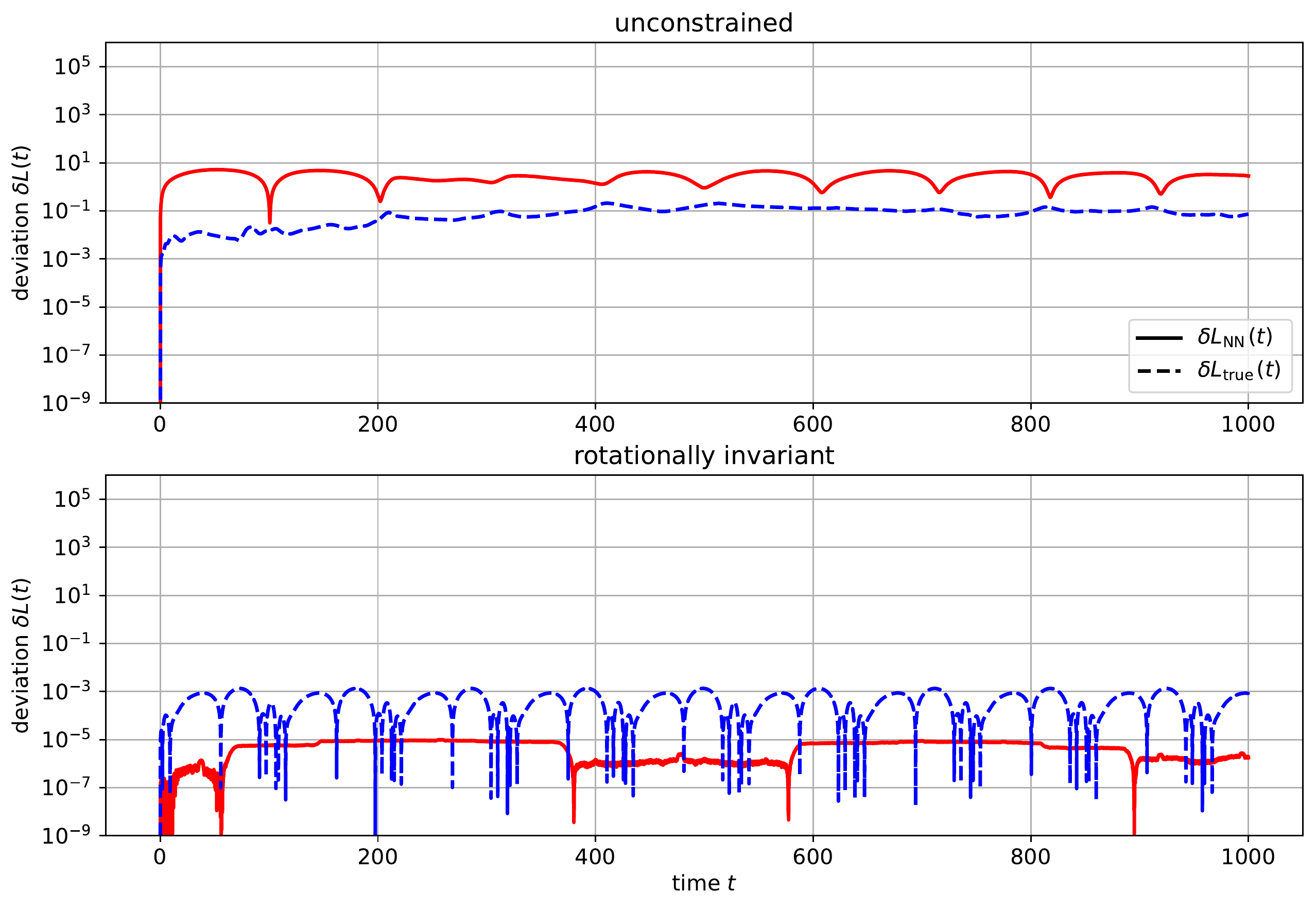}
    \caption{Conservation of angular momentum for the motion of a massive relativistic particle in the Schwarzschild metric as defined in Section \ref{sec:schwarzschild_problem} without (top) and with (bottom) constraints on the neural network to enforce rotational invariance of the Lagrangian. In each case the deviation $\delta L_{\text{NN}}(t)$ and $\delta L_{\text{true}}(t)$ of the absolute value of the total three dimensional angular momentum vectors are shown as a function of time. Both the ``neural network'' angular momentum $L_{\text{NN}}(t)$ defined in Eq. \eqref{eqn:relativistic_angular_momentum_NN} and the true angular momentum $L_{\text{true}}(t)$ defined in Eq. \eqref{eqn:relativistic_angular_momentum_true} are considered.}
    \label{fig:conservation_schwarzschild}
  \end{center}
\end{figure}

\begin{figure}
  \begin{center}
    \includegraphics[width=\linewidth]{\figdir/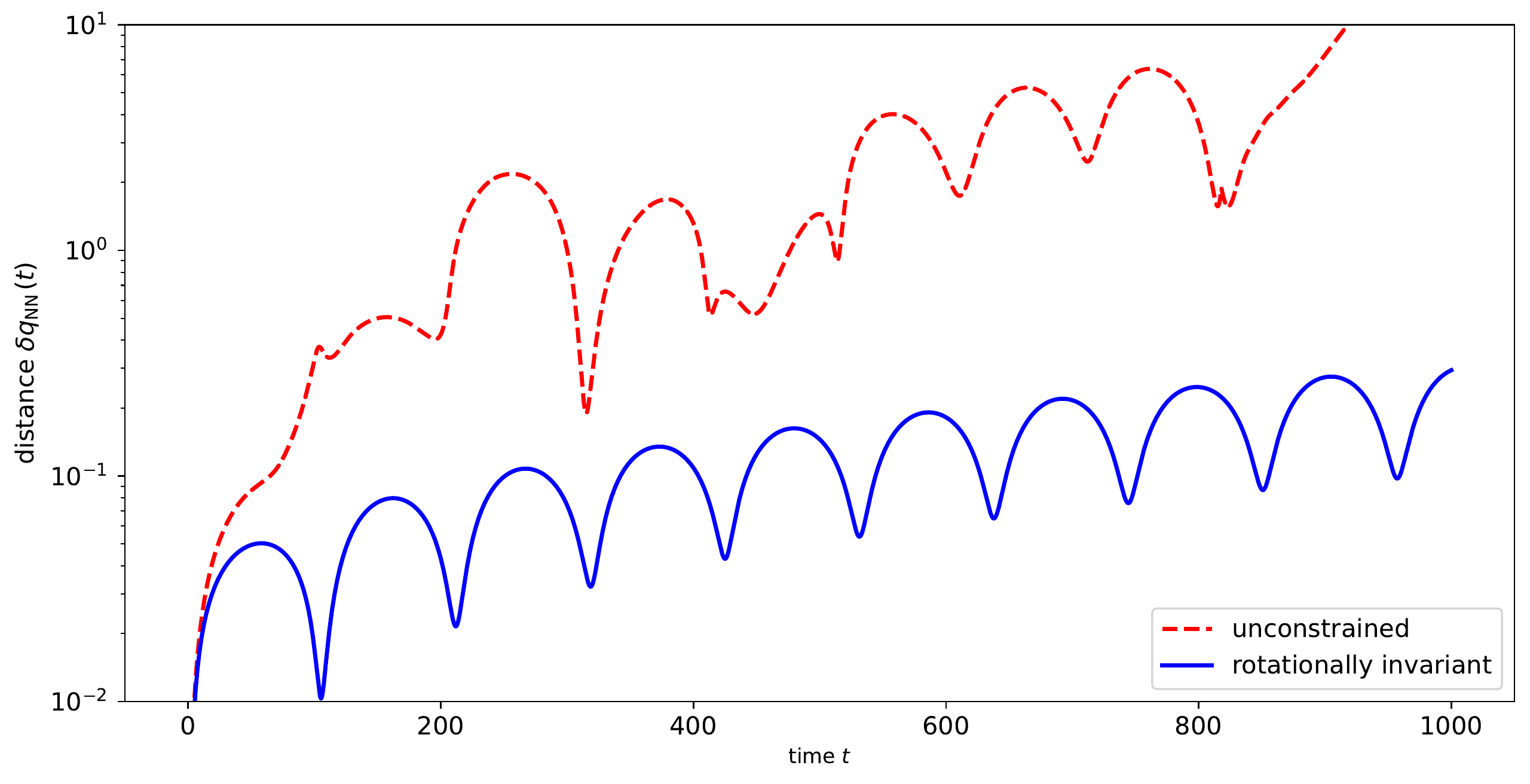}
    \caption{Evolution of the distance $\delta q_{\text{NN}}(t)=||q_{\text{NN}}(t)-q_{\text{NN}}^{\text{(perturbed)}}(t)||_2$ between two trajectories obtained with slightly different initial conditions for the motion of a massive relativistic particle in the Schwarzschild metric. The initial conditions that differ by $\delta q_{\text{NN}}(0)\sim 10^{-3}$.}
    \label{fig:pertubation_schwarzschild}
  \end{center}
\end{figure}

Fig. \ref{fig:trajectories_schwarzschild} shows the trajectories predicted with the trained neural network Lagrangian up to the final time $T=1000$. Again, a RK4 integrator with a timestep size $\Delta t=10^{-2}$ is used. The trajectory obtained with a slightly perturbed initial condition is also shown as a dashed line. In contrast to the non-relativistic equivalent, Fig. \ref{fig:trajectories_schwarzschild} (top) shows that for the unconstrained network the trajectories do not become unstable. However, they still diverge strongly from the true solution (shown as a dotted red line) and they oscillate about the $x-y$ plane, which indicates that angular momentum is not conserved. Perturbing the initial conditions also leads to very different trajectories at later times. Fig. \ref{fig:trajectories_schwarzschild} (bottom) demonstrates that the picture is fundamentally different for the rotationally invariant Lagrangian neural network. Here the trajectory generated with the neural network stays very close to the true solution (which is indeed almost completely hidden beneath the solid blue curve) and the two trajectories with slightly differing initial conditions only diverge slowly at later times. Furthermore, the motion appears to be completely confined to the $x-y$ plane, which indicates that the two corresponding components of the angular momentum are conserved.
Again, we investigate the conservation of angular momentum quantitatively by plotting the time evolution of the two-norms $\delta L_{\text{NN}}(t)$ and $\delta L_{\text{true}}(t)$ in Eq. \eqref{eqn:deltaL_definition}, now using the expressions for $L_{\text{NN}}$ and $L_{\text{true}}$ defined via Eqs. \eqref{eqn:relativistic_angular_momentum_true} and \eqref{eqn:relativistic_angular_momentum_NN}. As Fig. \ref{fig:conservation_schwarzschild} shows, the unconstrained network is able to learn the conservation of angular momentum to some degree, but the true angular momentum $L_{\text{true}}(t)$ deviates from its initial value by around $10\%$ at later times. This is consistent with the fact that the trajectories in Fig. \ref{fig:trajectories_schwarzschild} (top) are visibly not constrained to the $x-y$ plane, but do not stray too far from it either. The neural network with built-in rotational invariance, on the other hand, is able to reduce the relative deviation $\delta L_{\text{true}}(t)$ from the initial angular momentum to less than $10^{-3}$. As expected, the quantity $L_{\text{NN}}$, which would be zero in exact arithmetic and for an exact time integrator, is very small and never exceeds a value of around $10^{-5}$, which is consistent with (accumulated) rounding errors in single precision arithmetic.
Fig. \ref{fig:pertubation_schwarzschild} shows the distance $\delta q_{\text{NN}}(t) = ||q_{\text{NN}}(t)-q_{\text{NN}}^{\text{(perturbed)}}(t)||_2$ between the two neural network trajectories with slightly perturbed initial conditions. The plot confirms that the rotationally invariant Lagrangian is much more robust under perturbations of the initial conditions. Compared to the unconstrained Lagrangian neural network, the distance between the perturbed and the unperturbed trajectories is around one order of magnitude smaller and grows only moderately over time.
\subsubsection{Two interacting particles in $D=4$ dimensions}
\begin{figure}
  \begin{center}
    \includegraphics[width=\linewidth]{\figdir/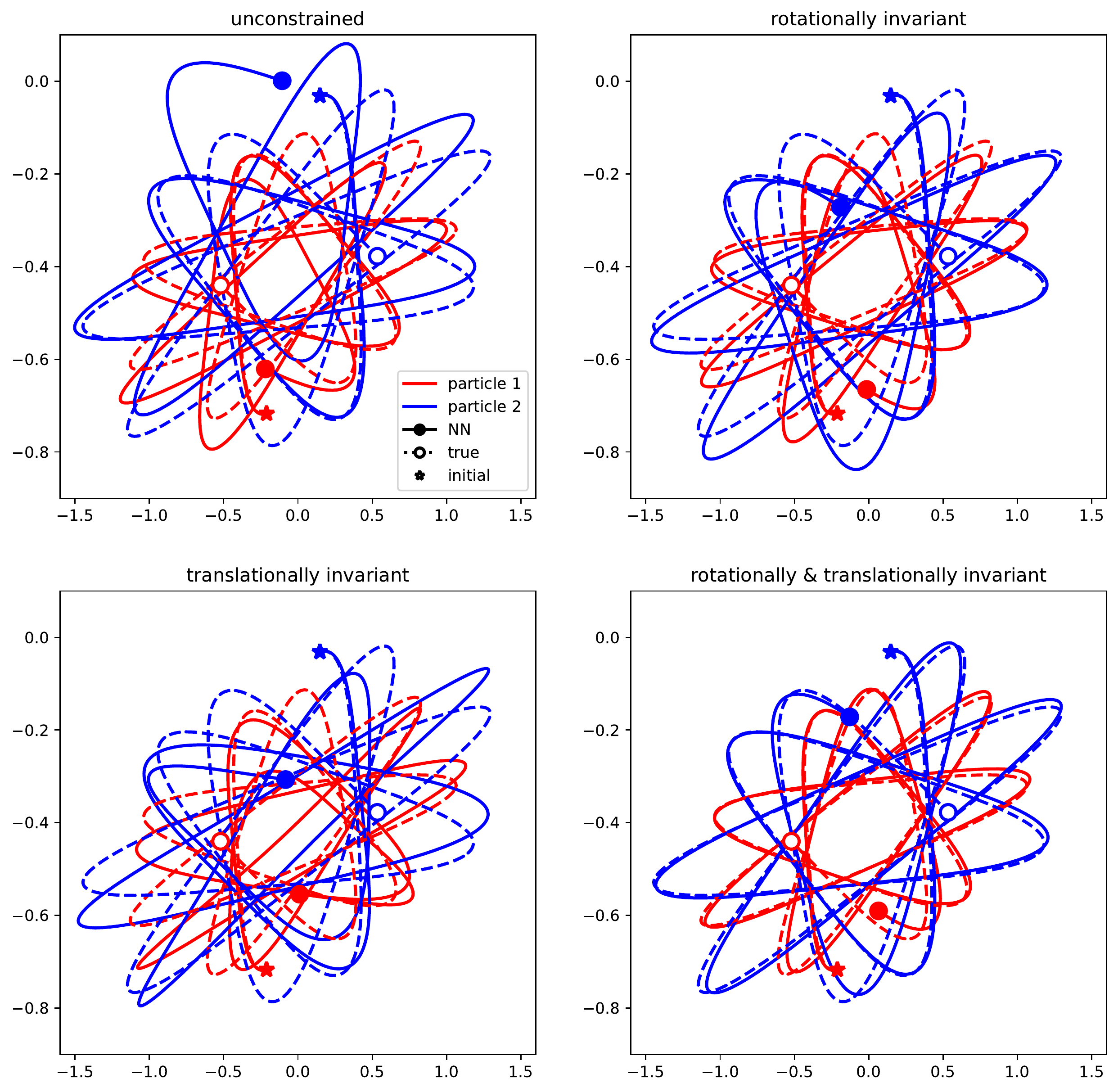}\\
    \caption{Trajectories for the two-particle problem assuming different constraints on the neural network Lagrangian as laid out in Tab. \ref{tab:lagrangian_restrictions}. In each case, both the true trajectories (dashed lines) and those obtained with the neural network Lagrangian (solid lines) are shown for the two particles. The starting positions are denoted with a star ($\medblackstar$) and the final positions are marked by circles ($\medblackcircle$/$\medcircle$). The figures shows a projection on the first two coordinates of four dimensional space.}
    \label{fig:trajectories_two_particle}
  \end{center}
\end{figure}

\begin{figure}
  \begin{center}
    \includegraphics[width=\linewidth]{\figdir/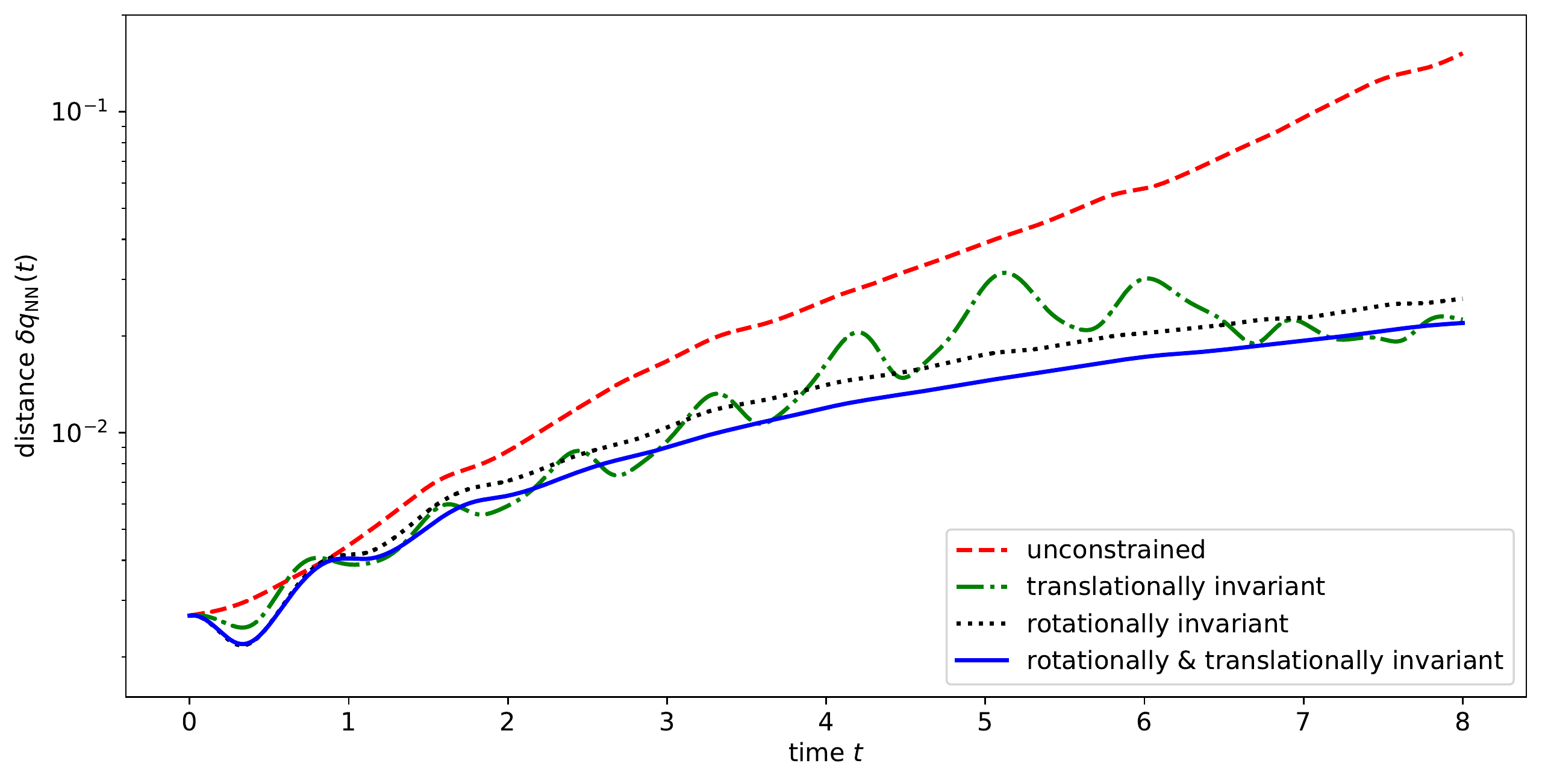}
    \caption{Evolution of the distance $\delta q_{\text{NN}}(t)= ||q_{\text{NN}}(t)-q_{\text{NN}}^{\text{(perturbed)}}(t)||_2$ between two trajectories obtained with slightly different initial conditions for the two-particle problem. The initial conditions that differ by $\delta q_{\text{NN}}(0)\sim 10^{-3}$.}
    \label{fig:pertubation_two_particle}
  \end{center}
\end{figure}

Finally, we consider the system of two interacting particles described in Section \ref{sec:two_particle_problem} for $D=4$ dimensions. The masses of the two particles are set to $m_1=1$ and $m_2=0.8$, while the parameters of the double well potential in Eq. \eqref{eqn:double_well} are fixed to $\mu=\kappa=1$. The initial condition was chosen such that at time $t=0$ the total linear momentum $M=m_1x^{(1)}+m_2x^{(2)}$ is zero. We consider four setups, corresponding to the constraints on the Lagrangian listed in the four quadrants of Tab. \ref{tab:lagrangian_restrictions}:
\begin{enumerate}
  \item Unconstrained Lagrangian
  \item Rotationally invariant Lagrangian
  \item Translationally invariant Lagrangian
  \item Rotationally and translationally invariant Lagrangian
\end{enumerate}
Since (except in the first case) the input to the neural network is first passed through a symmetry-enforcing layer, the number of inputs to the first hidden layer $D_1$ depends on which symmetries we assume for the Lagrangian, as shown in Tab. \ref{tab:n_inputs}. Note that as we do not enforce invariance under reflections (i.e. we only consider the $SO(4)$ subgroup instead of the full $O(4)$ group), in the case of a rotationally (but not necessarily translationally) invariant Lagrangian, the contraction $\varepsilon^{\alpha\beta\rho\sigma}x_{\alpha}^{(1)}x_{\beta}^{(2)}\dot{x}_{\rho}^{(1)}\dot{x}_{\sigma}^{(2)}$ with the Levi-Civita symbol is included in the set $\mathscr{R}\left(\{x^{(1)},x^{(2)},\dot{x}^{(2)},\dot{x}^{(2)}\}\right)$.
\begin{table}
  \begin{center}
    \begin{tabular}{|cc|>{\centering\arraybackslash}p{12ex}|>{\centering\arraybackslash}p{12ex}|}
      \hline
                      &     & \multicolumn{2}{c|}{rotationally invariant?}       \\
                      &     & no                                           & yes \\
      \hline
      translationally & no  & 16                                           & 11  \\
      \cline{2-4}
      invariant?      & yes & 12                                           & 6   \\
      \hline
    \end{tabular}
    \caption{Number of scalar inputs to the first dense layer $D_1$ of the Lagrangian neural network.}
    \label{tab:n_inputs}
  \end{center}
\end{table}
As in the other two cases, we generate synthetic training samples $(X^{(j)},\widehat{y}^{(j)})\in \mathbb{R}^{16}\times \mathbb{R}^{8}$ by integrating the true equations of motion obtained from the Lagrangian in Eq. \eqref{eqn:true_lagrangian_two_particle} with a RK4 integrator ($\Delta t=10^{-2}$) and adding normally distributed noise with a standard deviation of $\sigma=10^{-3}$.
Fig. \ref{fig:trajectories_two_particle} shows a projection of the four-dimensional trajectories onto the first two dimensions for each of the four considered setups. The true solution (dashed) and the neural network solution (solid) are shown for both particles with the final positions at time $T=8$ marked by circles ($\medblackcircle$/$\medcircle$). Visually, the completely unconstrained network gives the worst solution. While enforcing rotational or translation invariance improves this somewhat, the best qualitative agreement is achieved with the neural network Lagrangian that is both rotationally and translationally invariant. Although even in this case the final positions have a distance of order 1, this appears to be mainly attributed to the fact that the neural network solution lags behind the true solution. Ignoring this phase error, the trajectories generated with the neural network Lagrangian show good agreement with the true solution.
To quantify the conservation of the linear- and angular momentum, we again plot the time evolution of $\delta L_{\text{NN}}(t)$ and $\delta L_{\text{true}}$ as defined in Eq. \eqref{eqn:deltaL_definition}, but now with $L_{\text{NN}}$ and $L_{\text{true}}$ being the six-dimensional vectors defined in Eq. \eqref{eqn:two_particle_momentum_vectors}. In analogy to Eq. \eqref{eqn:deltaL_definition}, we define the time dependent deviation of the linear momentum as
\begin{xalignat}{2}
  \delta M_{\text{NN}}(t) & =\frac{\left|\left|M_{\text{NN}}(t)-M_{\text{NN}}(0)\right|\right|_2}{\left|\left||\widetilde{M}_{\text{NN}}(0)\right|\right|_2}, &
  \delta M_{\text{true}}(t) & =\frac{\left|\left|M_{\text{true}}(t)-M_{\text{true}}(0)\right|\right|_2}{\left|\left|\widetilde{M}_{\text{true}}(0)\right|\right|_2}.
\end{xalignat}
where $M_{\text{NN}}$ and $M_{\text{true}}$ are the four-dimensional linear momentum vectors given in Eq. \eqref{eqn:two_particle_momentum_vectors}. Since the initial conditions are chosen such that the total linear momentum is zero, we cannot normalise by $\left|\left|M_{\text{NN}}(0)\right|\right|_2$ and $\left|\left|M_{\text{true}}(0)\right|\right|_2$ but instead divide by the two-norm of the linear momentum of the first particle defined as
\begin{xalignat}{2}
  \widetilde{M}_{\text{NN}} &= \frac{\partial \lag_{\text{NN}}}{\partial \dot{x}^{(1)}}, &
  \widetilde{M}_{\text{true}} &= m^{(1)}\dot{x}^{(1)}.
\end{xalignat}
While the unconstrained network is able to learn the conservation of $L_{\text{true}}$ and $M_{\text{true}}$ to some degree, at the final time both quantities have deviated from their initial values by $10\%-100\%$. The picture is somewhat better once either translational or rotational invariance is enforced. In particular rotational invariance limits both the relative deviation of both the linear and angular momentum to below $\approx 10\%$. As expected, the linear and angular momenta $M_{\text{NN}}$ and $L_{\text{NN}}$ derived from the neural network Lagrangian (see Eq. \eqref{eqn:two_particle_NN_momenta}) are conserved up to rounding errors caused by single precision arithmetic. Enforcing both rotational and translation invariance improves the conservation of the true linear and angular momentum substantially, with both now not deviating by more than $1\%$ from their initial values. In contrast to the other cases, there also does not appear to be any upwards trend for $\delta M_{\text{true}}(t)$ or $\delta L_{\text{NN}}(t)$ as time increases.
\begin{figure}
  \begin{center}
    \includegraphics[width=\linewidth]{\figdir/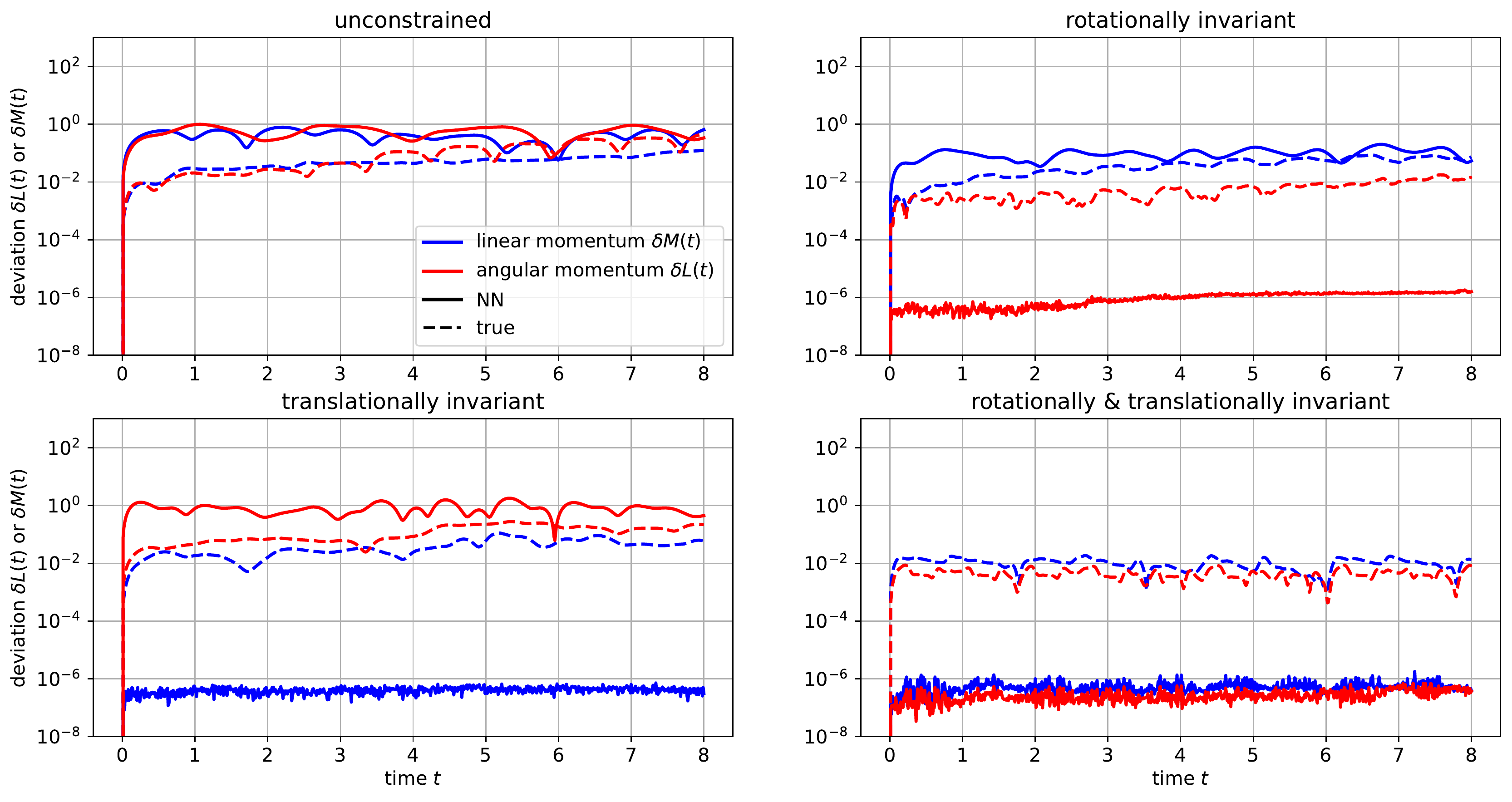}
    \caption{Conservation of linear and angular momentum for the four-dimensional two particle problem described in Section \ref{sec:two_particle_problem} with different constraints on the neural network to enforce rotational and/or translational invariance of the Lagrangian. In each case the relative deviations $\delta L_{\text{NN}}(t)$, $\delta M_{\text{NN}}(t)$, $\delta L_{\text{true}}(t)$ and $\delta M_{\text{true}}(t)$  of the linear and angular momentum vectors in Eq. \eqref{eqn:two_particle_momentum_vectors} are shown as a function of time. Both the ``neural network'' momenta defined in Eq. \eqref{eqn:two_particle_NN_momenta} and the true momenta in Eq. \eqref{eqn:momentum_true_twoparticle} are considered.}
    \label{fig:conservation_two_particle}
  \end{center}
\end{figure}
Finally, we study the stability of the neural network trajectories under small perturbations of the initial conditions. For this, in each of the four cases the initial position and velocity were perturbed by the same normal noise with standard deviation $10^{-3}$. The distance $\delta q_{\text{NN}}(t) = ||q_{\text{NN}}(t)-q_{\text{NN}}^{(\text{perturbed})}(t)||_2$ between the trajectories with these two different initial conditions is visualised in Fig. \ref{fig:pertubation_two_particle}. As can be seen from this plot, enforcing rotational and translational invariance reduces the growth of the perturbation significantly. However, only enforcing either rotational or translational invariance seems to be almost as good, with the latter performing even slightly better at large times $t$.
\section{Conclusion}\label{sec:conclusion}
In this paper we showed how the Lagrangian neural network approach in \cite{Cranmer2020} can be combined with Noether's Theorem to enforce conservation laws while learning the dynamics of a physical system from noisy data. As the numerical results demonstrate, this generates significantly more realistic trajectories since the network does not have to learn the underlying conservation laws from data. As expected, for each symmetry of the Lagrangian the neural network integrators conserve a quantity up to rounding errors from single precision arithmetic. While the corresponding invariant quantity of the true Lagrangian is not conserved exactly, its deviations from the initial value are significantly reduced. Using a Lagrangian with built-in symmetries also improves the stability of the solution in the sense that the system is less sensitive to small perturbations of the initial conditions.

There are several avenues for extending the work presented here. To demonstrate the principle ideas behind our approach we assumed that essentially unlimited (synthetic) data with relatively small errors is available for training. It would be interesting to instead train the Lagrangian neural networks on measured experimental data for real physical systems. While we already assumed that the training data is perturbed by random noise, in real-life situations measured data might be corrupted or only projections of the state space vectors might be available. For example, in celestial dynamics usually only motion perpendicular to the line of sight of the observer can be measured. More generally, it would be interesting to further explore the how the symmetry enforcing layers improve the solution for limited data with larger errors; initially one could simply vary the strength of the parameter $\sigma$ that characterises the strength of the noise in our numerical experiments.

So far, we considered Lagrangians that describe the motion of point particles. However, the Lagrangian formalism can also be extended to continuous systems such as (classical) field theories and it would be interesting to see how our approach works in this case. For example, in \cite{Cranmer2020} a discretised one-dimensional wave equation is studied and the authors show that the Lagrangian approach leads to the approximate conservation of energy.

It is not obvious that the methods described in this paper will still be practical when applied to larger systems of many interacting particles. One concern here is that the need to take derivatives of the Lagrangian will amplify errors. Since the Lagrangian depends on the positions and velocities of all particles, this issue might become more severe as the size of the system grows. In addition, when simulating $N$ particles in $d$ dimensions the $Nd\times Nd$ Hessian matrix $J_{\dot{q},\dot{q}}$ in Eq. \eqref{eqn:J_matrices} needs to be inverted, which is expensive (the cost of the inversion grows as $\mathcal{O}((Nd)^3)$) and might not be numerically stable for large $N$. The Hamiltonian formulation avoids this issue but, as argued in \cite{Cranmer2020}, it requires training data in canonical coordinates. If the Hamiltonian is unknown (consider a real-world situation in which we want to infer the unknown physical laws from measurements), it might not be obvious what these canonical coordinates are in the first place. If the most general invariant Lagrangian is to be constructed, the number of features after the symmetry enforcing layer will grow rapidly as $\mathcal{O}(N^2)$; see, however, \cite[Section 5]{Villar2021} which discusses a bound on the number of independent scalars.

In addition, only assuming invariance under continuous symmetries might not be sufficient to learn the dynamics for limited data. Building additional assumptions into the neural network architecture will help with this. To reduce the number of features (= scalar invariant returned by the symmetry enforcing layer), one could introduce additional constraints, such as the assumption that the Lagrangian is invariant under permutations of the particles \cite{Villar2021} (provided they have the same mass), the interactions are local or that the kinetic energy is of the form
\begin{equation}
  T(\dot{x}^{(1)},\dots,\dot{x}^{(N)}) = \sum_{j=1}^N \frac{1}{2}m_j \left(\dot{x}^{(j)}\right)^2\label{eqn:T_restricted}
\end{equation}
and the potential only depends on the coordinates (compare to \cite{Finzi2020}). This results in a constant, diagonal Hessian but it seriously restricts the class of systems that can be simulated: there are simple systems such as the double pendulum for which the Lagrangian is non-separable and the kinetic energy cannot be written in the form in Eq. \eqref{eqn:T_restricted}.

The focus of the work presented here has been on engineering features in the first hidden neural network layer, such that exact invariance of the Lagrangian under a continuous symmetry group is guaranteed. Since the systems for which we carried out numerical experiments are small, we did not focus on optimising the structure of subsequent layers. As discussed in the introduction, GNNs \cite{Battaglia2018,Battaglia2016} allow the successful simulation of large interacting particle systems and have been used to model the Hamiltonian in \cite{Sanchez2019}. It would therefore be very interesting to explore whether our approach can be combined with the GNN architecture. Possibly this could be achieved by using suitably engineered invariant features such as the squared distances $(x^{(j)}-x^{(k)})^2$ in the embeddings for edges connecting two particles $j$ and $k$ and $(\dot{x}^{(j)})^2$, $x^{(j)}\cdot \dot{x}^{(j)}$ etc. in the vertex embeddings, while returning the Hamiltonian or Lagrangian in the decoder as in \cite{Sanchez2019}. Since this introduces additional inductive biases, it might help with the simulation of significantly larger particle systems.

More generally, it will be interesting to compare the methods introduced in the present paper to other approaches such as \cite{Sanchez2020,Finzi2020,Sanchez2019,Schuett2017}.
\section*{Acknowledgements}
The author would like to thank Patrick Lavelle and Tony Shardlow for useful discussions during the early stages of this project. We are also grateful to George Em Karniadakis, Patrick Obin Sturm and the two anonymous reviewers for feedback that helped to put our work into a wider context.
\FloatBarrier
\appendix
\section{Rotationally invariant combinations of $\boldsymbol{n}$ vectors}\label{sec:rotationally_invariant_combinations}
Let $A=\{a^{(1)},a^{(2)},\dots,a^{(n)}\}$ with $a^{(\alpha)}\in\mathbb{R}^D$ be a set of $n$ vectors in $D$ dimensions. We want to construct the set $\mathscr{R}(A)$ of all rotationally invariant combinations of these $n$ vectors. For this observe that the only quantities that are invariant under the special orthogonal group $SO(D)$ are scalars which are obtained by either contracting pairs of vectors with the metric tensor $\eta$ or by contracting $D$-tuples of the vectors with the $D$-dimensional Levi-Civita symbol $\varepsilon$. The tensors $\eta$ and $\varepsilon$ are defined as follows:
\begin{equation}
  \begin{aligned}
    \eta^{ij}                     & = \begin{cases}
      1 & \text{if $i=j$}  \\
      0 & \text{otherwise}
    \end{cases}\qquad\text{for $i,j=1,2,\dots,D$,} \\
    \varepsilon^{i_1i_2\dots i_D} & = \begin{cases}
      \text{sgn}(\sigma) & \text{if $i_1,i_2,\dots,i_p$ is obtained from the numbers $1,2,\dots,D$ through the permutation $\sigma$} \\
      0                  & \text{otherwise,}
    \end{cases}
  \end{aligned}
  \label{eqn:levi_civita_definition}
\end{equation}
where $\text{sgn}(\sigma)\in\{-1,+1\}$ denotes the sign of the permutation $\sigma$.
There are ${n\choose 2}+n = \frac{1}{2}n(n+1)$ different contractions with the metric tensor, namely the dot products
\begin{equation}
  a^{(\alpha)}\cdot a^{(\beta)} := \eta^{ij}a^{(\alpha)}_ia^{(\beta)}_j
  \qquad\text{for $1\le \alpha\le \beta\le n$}.\label{eqn:metric_scalar_product}
\end{equation}
To simplify notation, we have adopted the Einstein sum convention of summing over pairs of identical upper and lower indices.
The number of full contractions with the Levi-Civita symbol is ${n\choose D}$. Each distinct full contraction is given by a choice of indices $\alpha_1,\alpha_2,\dots,\alpha_D$ with
$1\le\alpha_1<\alpha_2<\dots<\alpha_D\le D$ and it can be written down explicitly as
\begin{equation}
  a^{(\alpha_1)}:a^{(\alpha_2)}:\dots:a^{(\alpha_D)} := \varepsilon^{i_1i_2\dots i_D}
  a^{(\alpha_1)}_{i_1}a^{(\alpha_2)}_{i_2}\dots a^{(\alpha_D)}_{i_D}.\label{eqn:levi_civita_contractions}
\end{equation}
There are no full contractions with the Levi-Civita symbol if $n<D$ and there is only one full contraction for $n=D$, which is the special case considered in \cite[Theorem 2.9.A]{Weyl1946}. At this point it is also worth pointing out that here we only assume invariance under the \textit{special} orthogonal group $SO(D)\subset O(D)$. In addition to rotations, the orthogonal group $O(D)$ also contains reflections, which change the sign of the contractions in \eqref{eqn:levi_civita_contractions}. Hence, if we were to enforce this additional constraint, only even powers of these terms can appear in the Lagrangian. Since the product $\varepsilon^{i_1i_2\dots i_D}\varepsilon^{j_1j_2\dots j_D}$ of two Levi-Civita symbols can be expressed in terms of the metric tensor (see e.g. \cite[Chapter 6]{Cvitanovic2008}), any expressions of this form can then be reduced to scalar products, and it would be sufficient to only include terms of the form given in Eq. \eqref{eqn:metric_scalar_product} in $\mathscr{R}(A)$; this is also discussed in \cite[Section 3, Lemma 1]{Villar2021}. However, in this work we are only interested in continuous transformations which are the same in $SO(D)$ and $O(D)$ and we therefore include expressions of the form in Eq. \eqref{eqn:levi_civita_contractions} in $\mathscr{R}(A)$.
\section{Loss histories}\label{sec:loss_histories}
The following figures show the loss histories for the MSE error defined in \eqref{eqn:MSE_loss} for all three problems considered in this paper.
\begin{figure}
  \begin{center}
    \includegraphics[width=\linewidth]{\figdir/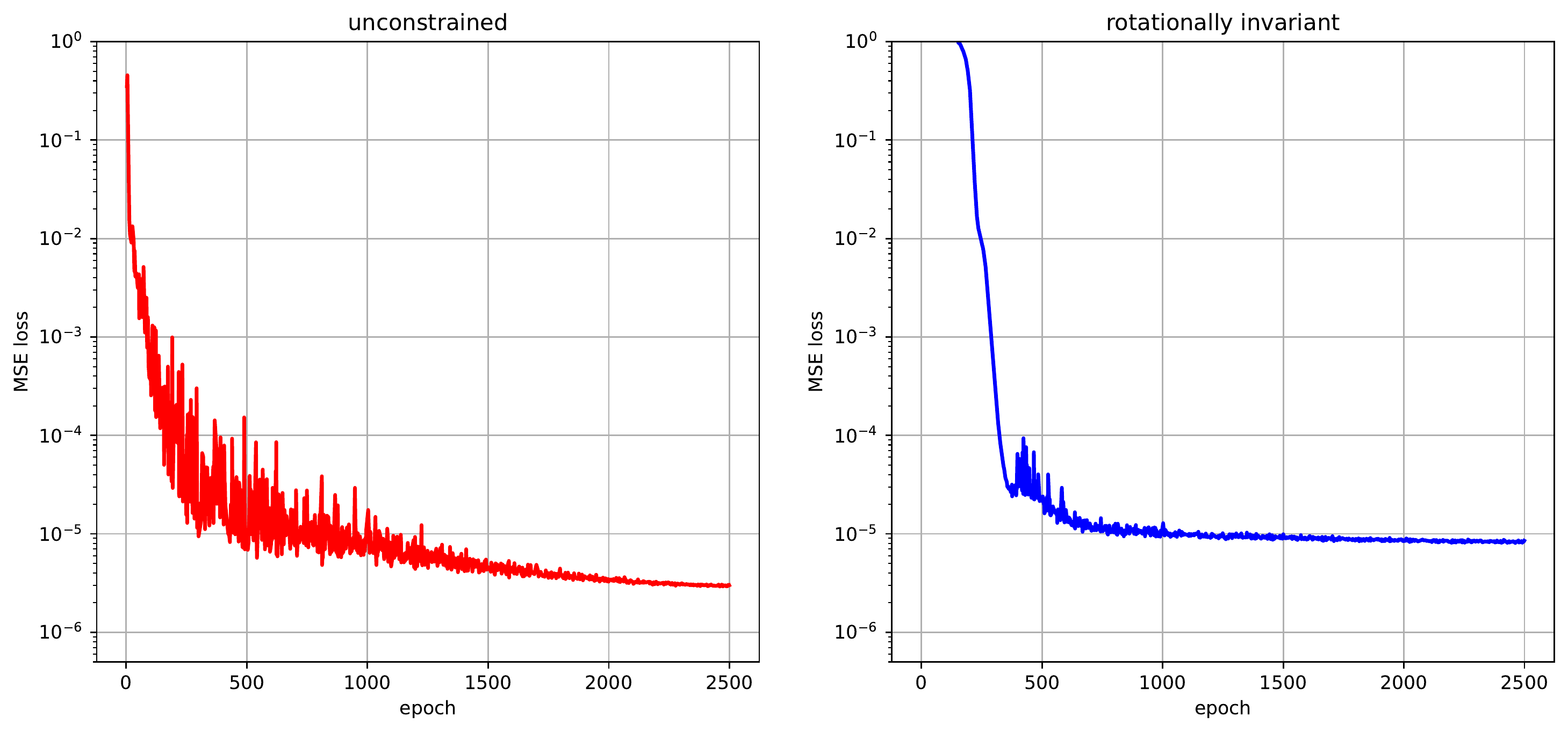}
    \caption{Evolution of the mean squared loss for the unconstrained (left) and the rotationally invariant (right) motion of a non-relativistic particle in a gravitational potential as defined in Section \ref{sec:kepler_problem}.}
    \label{fig:loss_history_kepler}
  \end{center}
\end{figure}
\begin{figure}
  \begin{center}
    \includegraphics[width=\linewidth]{\figdir/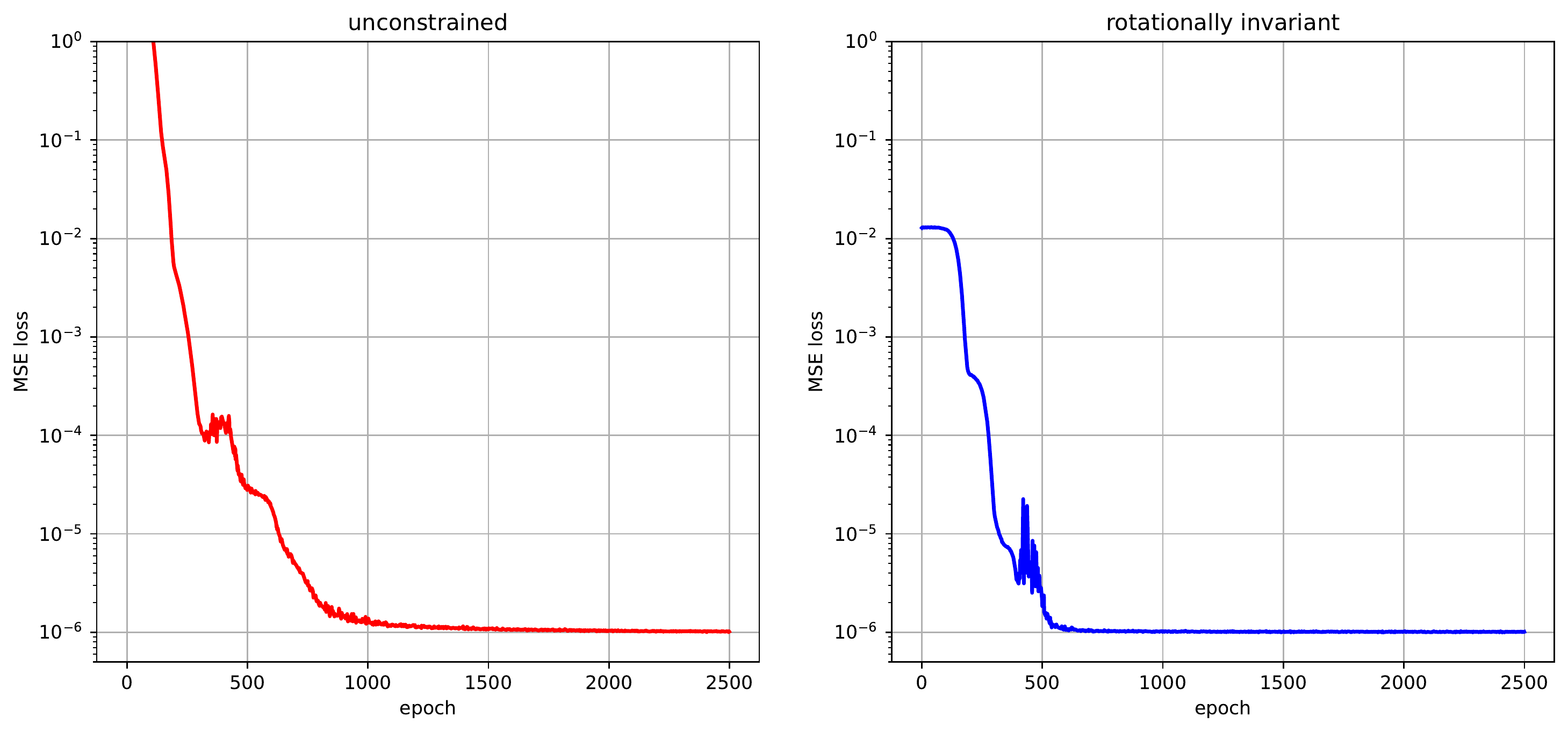}
    \caption{Evolution of the mean squared loss for the unconstrained (left) and the rotationally invariant (right) motion of a massive relativistic particle in the Schwarzschild metric as defined in Section \ref{sec:schwarzschild_problem}.}
    \label{fig:loss_history_schwarzschild}
  \end{center}
\end{figure}
\begin{figure}
  \begin{center}
    \includegraphics[width=\linewidth]{\figdir/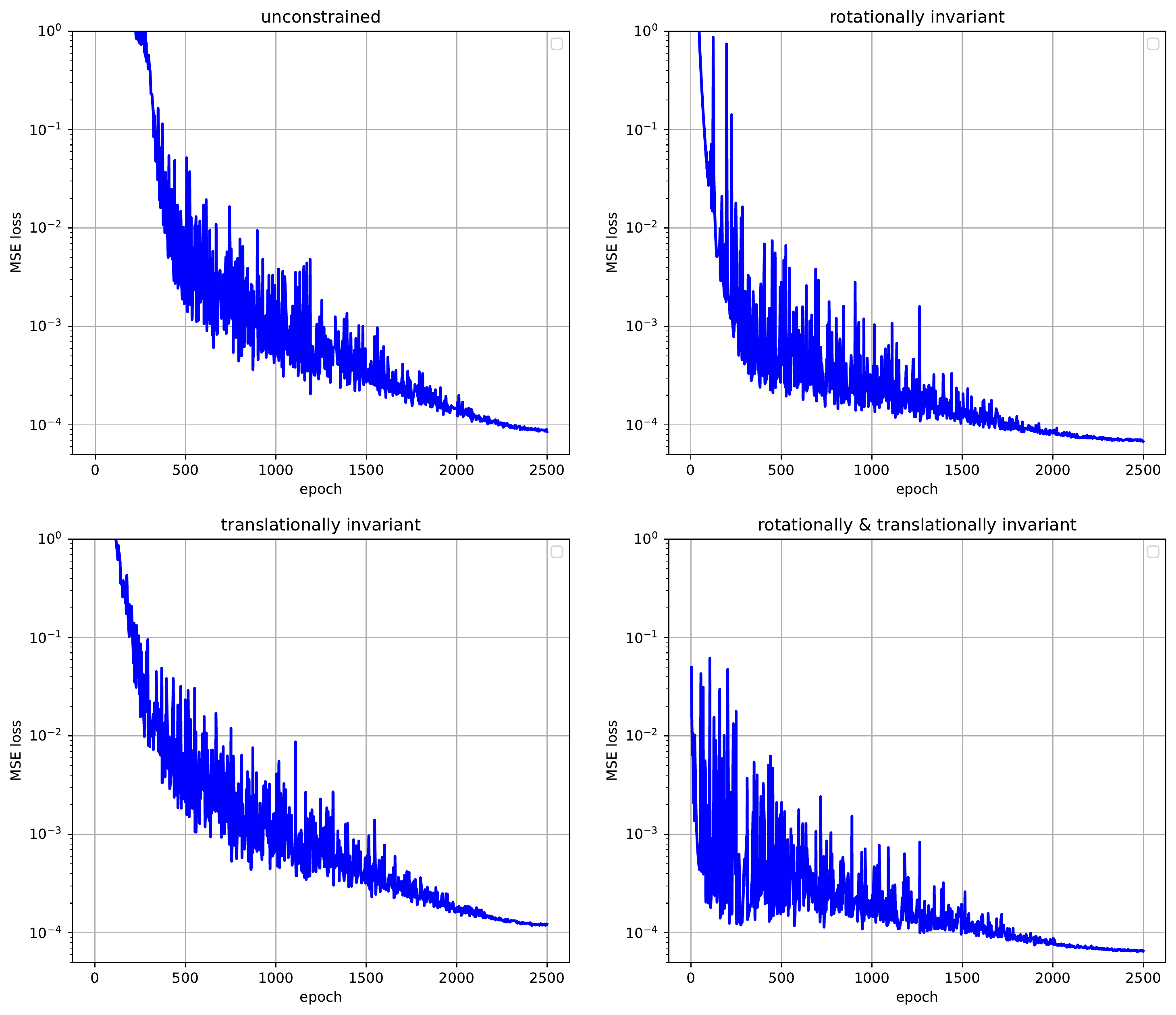}
    \caption{Evolution of the mean squared loss for the two particle problem defined in Section \ref{sec:two_particle_problem}, making different assumptions on the invariance of the neural network Lagrangian.}
    \label{fig:loss_history_two_particle}
  \end{center}
\end{figure}
\bibliographystyle{elsarticle-num}

\begin{thebibliography}{10}
  \expandafter\ifx\csname url\endcsname\relax
    \def\url#1{\texttt{#1}}\fi
  \expandafter\ifx\csname urlprefix\endcsname\relax\def\urlprefix{URL }\fi
  \expandafter\ifx\csname href\endcsname\relax
    \def\href#1#2{#2} \def\path#1{#1}\fi

  \bibitem{He2016}
  K.~He, X.~Zhang, S.~Ren, J.~Sun, {Deep residual learning for image
      recognition}, in: Proceedings of the IEEE conference on computer vision and
  pattern recognition, 2016, pp. 770--778.

  \bibitem{Ren2015}
  S.~Ren, K.~He, R.~Girshick, J.~Sun, {Faster R-CNN: Towards real-time object
      detection with region proposal networks}, Advances in neural information
  processing systems 28 (2015).

  \bibitem{Krizhevsky2012}
  A.~Krizhevsky, I.~Sutskever, G.~E. Hinton, {Imagenet classification with deep
      convolutional neural networks}, Advances in neural information processing
  systems 25 (2012).

  \bibitem{Einstein1923}
  A.~Einstein, {Die Grundlage der allgemeinen Relativit{\"a}tstheorie}, in: Das
  Relativit{\"a}tsprinzip, Springer, 1923, pp. 81--124.

  \bibitem{Gasser1984}
  J.~Gasser, H.~Leutwyler, {Chiral perturbation theory to one loop}, Annals of
  Physics 158~(1) (1984) 142--210.

  \bibitem{Scherer2003}
  S.~Scherer, {Introduction to chiral perturbation theory}, in: Advances in
  Nuclear Physics, Volume 27, Springer, 2003, pp. 277--538.

  \bibitem{Noether1918}
  E.~Noether, \href{http://eudml.org/doc/59024}{{Invariante Variationsprobleme}},
  {Nachrichten von der Gesellschaft der Wissenschaften zu G\"{o}ttingen,
  Mathematisch-Physikalische Klasse} 1918 (1918) 235--257.
  \newline\urlprefix\url{http://eudml.org/doc/59024}

  \bibitem{Pfaff2020}
  T.~Pfaff, M.~Fortunato, A.~Sanchez-Gonzalez, P.~W. Battaglia, {Learning
      mesh-based simulation with graph networks}, arXiv preprint arXiv:2010.03409
  (2020).

  \bibitem{Sanchez2020}
  A.~Sanchez-Gonzalez, J.~Godwin, T.~Pfaff, R.~Ying, J.~Leskovec, P.~Battaglia,
  {Learning to simulate complex physics with graph networks}, in: International
  conference on machine learning, PMLR, 2020, pp. 8459--8468.

  \bibitem{Kadupitiya2022}
  J.~Kadupitiya, G.~C. Fox, V.~Jadhao, {Solving Newton’s equations of motion
      with large timesteps using recurrent neural networks based operators},
  Machine Learning: Science and Technology 3~(2) (2022) 025002.

  \bibitem{Hochreiter1997}
  S.~Hochreiter, J.~Schmidhuber, {Long short-term memory}, Neural computation
  9~(8) (1997) 1735--1780.

  \bibitem{Hairer1993}
  E.~Hairer, S.~P. N{\o}rsett, G.~Wanner, {Solving ordinary differential
      equations. 1, Nonstiff problems}, Springer-Vlg, 1993.

  \bibitem{Greydanus2019}
  S.~Greydanus, M.~Dzamba, J.~Yosinski, {Hamiltonian neural networks}, Advances
  in neural information processing systems 32 (2019).

  \bibitem{Chen2019}
  Z.~Chen, J.~Zhang, M.~Arjovsky, L.~Bottou, {Symplectic recurrent neural
      networks}, arXiv preprint arXiv:1909.13334 (2019).

  \bibitem{Verlet1967}
  L.~Verlet, {Computer ``experiments'' on classical fluids. I. Thermodynamical
      properties of Lennard-Jones molecules}, Physical review 159~(1) (1967) 98.

  \bibitem{Cranmer2020}
  M.~Cranmer, S.~Greydanus, S.~Hoyer, P.~Battaglia, D.~Spergel, S.~Ho,
  {Lagrangian neural networks}, arXiv preprint arXiv:2003.04630 (2020).

  \bibitem{Mattheakis2019}
  M.~Mattheakis, P.~Protopapas, D.~Sondak, M.~Di~Giovanni, E.~Kaxiras, {Physical
      symmetries embedded in neural networks}, arXiv preprint arXiv:1904.08991
  (2019).

  \bibitem{Ling2016}
  J.~Ling, R.~Jones, J.~Templeton, {Machine learning strategies for systems with
      invariance properties}, Journal of Computational Physics 318 (2016) 22--35.

  \bibitem{Beucler2021}
  T.~Beucler, M.~Pritchard, S.~Rasp, J.~Ott, P.~Baldi, P.~Gentine, Enforcing
  analytic constraints in neural networks emulating physical systems, Physical
  Review Letters 126~(9) (2021) 098302.

  \bibitem{Sturm2020}
  P.~O. Sturm, A.~S. Wexler, {A mass-and energy-conserving framework for using
      machine learning to speed computations: a photochemistry example},
  Geoscientific Model Development 13~(9) (2020) 4435--4442.

  \bibitem{Sturm2022}
  P.~O. Sturm, A.~S. Wexler, {Conservation laws in a neural network architecture:
      enforcing the atom balance of a Julia-based photochemical model (v0. 2.0)},
  Geoscientific Model Development 15~(8) (2022) 3417--3431.

  \bibitem{Villar2021}
  S.~Villar, D.~W. Hogg, K.~Storey-Fisher, W.~Yao, B.~Blum-Smith, {Scalars are
      universal: Equivariant machine learning, structured like classical physics},
  Advances in Neural Information Processing Systems 34 (2021) 28848--28863.

  \bibitem{Weyl1946}
  H.~Weyl, {The classical groups: their invariants and representations}, no.~1,
  Princeton university press, 1946.

  \bibitem{Yao2021}
  W.~Yao, K.~Storey-Fisher, D.~W. Hogg, S.~Villar, {A simple equivariant machine
      learning method for dynamics based on scalars}, arXiv preprint
  arXiv:2110.03761 (2021).

  \bibitem{Thomas2018}
  N.~Thomas, T.~Smidt, S.~Kearnes, L.~Yang, L.~Li, K.~Kohlhoff, P.~Riley, {Tensor
      field networks: Rotation-and translation-equivariant neural networks for 3d
      point clouds}, arXiv preprint arXiv:1802.08219 (2018).

  \bibitem{Cohen2016}
  T.~Cohen, M.~Welling, {Group equivariant convolutional networks}, in:
  International conference on machine learning, PMLR, 2016, pp. 2990--2999.

  \bibitem{Finzi2020}
  M.~Finzi, S.~Stanton, P.~Izmailov, A.~G. Wilson, {Generalizing convolutional
      neural networks for equivariance to lie groups on arbitrary continuous data},
  in: International Conference on Machine Learning, PMLR, 2020, pp. 3165--3176.

  \bibitem{Battaglia2018}
  P.~W. Battaglia, J.~B. Hamrick, V.~Bapst, A.~Sanchez-Gonzalez, V.~Zambaldi,
  M.~Malinowski, A.~Tacchetti, D.~Raposo, A.~Santoro, R.~Faulkner, et~al.,
  {Relational inductive biases, deep learning, and graph networks}, arXiv
  preprint arXiv:1806.01261 (2018).

  \bibitem{Battaglia2016}
  P.~Battaglia, R.~Pascanu, M.~Lai, D.~Jimenez~Rezende, et~al., {Interaction
      networks for learning about objects, relations and physics}, Advances in
  neural information processing systems 29 (2016).

  \bibitem{Sanchez2019}
  A.~Sanchez-Gonzalez, V.~Bapst, K.~Cranmer, P.~Battaglia, {Hamiltonian graph
      networks with ode integrators}, arXiv preprint arXiv:1909.12790 (2019).

  \bibitem{Zhang2022}
  Z.~Zhang, Y.~Shin, G.~Em~Karniadakis, {GFINNs: GENERIC formalism informed
      neural networks for deterministic and stochastic dynamical systems},
  Philosophical Transactions of the Royal Society A 380~(2229) (2022) 20210207.

  \bibitem{Grmela1997}
  M.~Grmela, H.~C. {\"O}ttinger, {Dynamics and thermodynamics of complex fluids.
      I. Development of a general formalism}, Physical Review E 56~(6) (1997) 6620.

  \bibitem{Oettinger1997}
  H.~C. {\"O}ttinger, M.~Grmela, {Dynamics and thermodynamics of complex fluids.
      II. Illustrations of a general formalism}, Physical Review E 56~(6) (1997)
  6633.

  \bibitem{Oettinger2005}
  H.~C. {\"O}ttinger, {Beyond equilibrium thermodynamics}, John Wiley \& Sons,
  2005.

  \bibitem{Oettinger2018}
  H.~C. {\"O}ttinger, {GENERIC integrators: structure preserving time integration
      for thermodynamic systems}, Journal of Non-Equilibrium Thermodynamics 43~(2)
  (2018) 89--100.

  \bibitem{Hernandez2021}
  Q.~Hern{\'a}ndez, A.~Bad{\'\i}as, D.~Gonz{\'a}lez, F.~Chinesta, E.~Cueto,
  {Structure-preserving neural networks}, Journal of Computational Physics 426
  (2021) 109950.

  \bibitem{Lee2021}
  K.~Lee, N.~Trask, P.~Stinis, {Machine learning structure preserving brackets
      for forecasting irreversible processes}, Advances in Neural Information
  Processing Systems 34 (2021) 5696--5707.

  \bibitem{Arnold2013}
  V.~I. Arnold, {Mathematical methods of classical mechanics}, Vol.~60, Springer
  Science \& Business Media, 2013.

  \bibitem{Schwarzschild1916}
  K.~Schwarzschild, {\"{U}ber das Gravitationsfeld eines Massenpunktes nach der
  Einstein'schen Theorie}, Berlin. Sitzungsberichte 18 (1916).

  \bibitem{Droste1917}
  J.~Droste, {The field of a single centre in Einstein's theory of gravitation,
      and the motion of a particle in that field}, Ned. Acad. Wet., SA 19 (1917)
  197.

  \bibitem{Tensorflow2015}
  M.~Abadi, A.~Agarwal, P.~Barham, E.~Brevdo, Z.~Chen, C.~Citro, G.~S. Corrado,
  A.~Davis, J.~Dean, M.~Devin, S.~Ghemawat, I.~Goodfellow, A.~Harp, G.~Irving,
  M.~Isard, Y.~Jia, R.~Jozefowicz, L.~Kaiser, M.~Kudlur, J.~Levenberg,
  D.~Man\'{e}, R.~Monga, S.~Moore, D.~Murray, C.~Olah, M.~Schuster, J.~Shlens,
  B.~Steiner, I.~Sutskever, K.~Talwar, P.~Tucker, V.~Vanhoucke, V.~Vasudevan,
  F.~Vi\'{e}gas, O.~Vinyals, P.~Warden, M.~Wattenberg, M.~Wicke, Y.~Yu,
  X.~Zheng, \href{https://www.tensorflow.org/}{{TensorFlow}: Large-scale
    machine learning on heterogeneous systems}, software available from
  tensorflow.org (2015).
  \newline\urlprefix\url{https://www.tensorflow.org/}

  \bibitem{code_release}
  E.~H. Mueller, \href{https://doi.org/10.5281/zenodo.7108069}{{Code for paper on
        "Exact conservation laws for neural network integrators of dynamical
        systems"}} (Sep. 2022).
  \newblock \href {https://doi.org/10.5281/zenodo.7108069}
  {\path{doi:10.5281/zenodo.7108069}}.
  \newline\urlprefix\url{https://doi.org/10.5281/zenodo.7108069}

  \bibitem{Schuett2017}
  K.~Sch{\"u}tt, P.-J. Kindermans, H.~E. Sauceda~Felix, S.~Chmiela,
  A.~Tkatchenko, K.-R. M{\"u}ller, {Schnet: A continuous-filter convolutional
      neural network for modeling quantum interactions}, Advances in neural
  information processing systems 30 (2017).

  \bibitem{Cvitanovic2008}
  P.~Cvitanovi{\'c}, {Group theory}, in: Group Theory, Princeton University
  Press, 2008.

\end{thebibliography}

\end{document}